\documentclass[11pt]{article}
\usepackage{theorem,hyperref}
\usepackage{graphicx}

\def \build#1#2#3{\mathrel{\mathop{\kern 0pt#1}\limits_{#2}^{#3}}}

\def \tend{\longrightarrow}
\def \be{\begin{eqnarray*}}
\def \ee{\end{eqnarray*}}
\def \B{\Big}
\def \pr{\indent \sl Proof : \rm}

\newcommand \bel{\begin{eqnarray}\label}
\newcommand \eel{\end{eqnarray}}

\newcommand\mathbb{}

\newtheorem{lem}{\indent Lemma}[section]

\newtheorem{pro}[lem]{\indent Proposition}
\newtheorem{cor}[lem]{\indent Corollary}
\newtheorem{theo}[lem]{\indent Theorem}

\def\rompar(#1){\textup(#1\textup)}    

\begin{document}

\begin{center}

\bf\LARGE{Phase transition
for parking blocks,\\ Brownian excursion
and coalescence\footnote{Random Structures Algorithms 21 (2002), no. 1, 76-119.}\\}
\rm
\vspace{2mm}
\Large
P. Chassaing \footnote{Institut Elie Cartan,
INRIA, CNRS and Universit\'e Henri Poincar\'e,
\\
BP 239,\ \ 54 506 Vandoeuvre  Cedex, France.
\\
chassain@iecn.u-nancy.fr}
\& G. Louchard
\footnote{Universit\'e Libre de Bruxelles,
D\'epartement d'Informatique,
\\Campus Plaine, CP 212,
  Bvd du Triomphe,
1050 Bruxelles, Belgium.\\
louchard@ulb.ac.be}\\

\vspace{3mm}

\rm
\vspace{3mm}
\rm

\vspace{2cm}

\end{center}

\begin{quotation}\small
\noindent
\textbf{Abstract. } In this paper, we consider
hashing with linear probing for a hashing table with $m$ places,
 $n$ items ($n<m$), and $\ell =m-n$ empty places. For a
 non computer science-minded
 reader, we shall use
 the metaphore of $n$ cars parking on $m$ places:
each car $c_i$ chooses a place $p_i$ at random, and if $p_i$
is occupied, $c_i$ tries successively $p_i+1$, $p_i+2$, until it finds
an empty place.
Pittel  \cite{PITT} proves that when $\ell/m$ goes
 to some positive limit $\beta<1$,
the size $B_1^{m,\ell}$ of the largest block of
consecutive cars satisfies
$2(\beta-1-\log\beta)B_1^{m,\ell}=2\log m-3
\log\log m+\Xi_m$,
 where $\Xi_m$ converges weakly to an extreme-value distribution.
In this paper we examine  at which level for $n$
a phase transition occurs  between  $B_1^{m,\ell}=o(m)$ and
$m-B_1^{m,\ell}=o(m)$. The intermediate case
 reveals an interesting
behaviour of sizes of blocks, related
 to the standard additive coalescent in the same way
as the sizes of connected components
of the random graph are related to the multiplicative coalescent.
\vspace{2mm}


\vspace{2mm}

\noindent
\textbf{A.M.S. Classification. }
60C05, 60J65, 60F05, 68P10, 68R05.
\vspace{3,5cm}

\end{quotation}

\vspace{1,8cm}

\section{Main results}
\label{mr}

\subsection{Emergence of a giant block}
\label{giant}

  We consider
hashing with linear probing for a hashing
 table with a set of $m$ places,
$\{1,2, ..., m\}$,
 $n$ items $\{c_1,c_2, ..., c_{n}\}$,
 and $\ell=m-n$ empty places ($\ell>0$). Hashing with linear probing
is a fundamental object
 in analysis of algorithms: its study goes back to the 1960's
\cite{KNU,KHW}
and is still active
\cite{AHASH,HASH,KNU2,PITT}. For a
 non computer science-minded
 reader, we shall use, all along the paper,
 the metaphore of $n$ cars parking on $m$ places,
leaving $\ell$ places empty:
each car $c_i$ chooses a place $p_i$ at random, and if $p_i$
is occupied, $c_i$ tries successively $p_i+1$, $p_i+2$,
until it finds
an empty place. We use the convention
  that place $m+1$ is also place $1$.

\begin{figure}[ht]
\hspace{116pt}
\includegraphics[width=6.7cm]{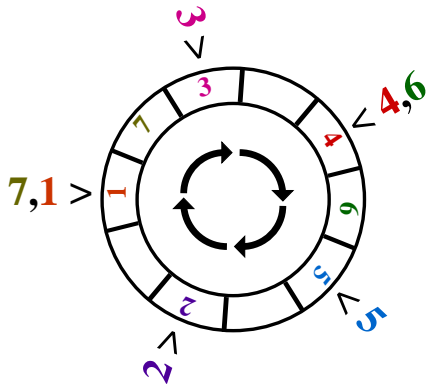}
\vspace{-1,3cm}
\begin{center}\caption{An example where  $(m,n)=(10,7)$ and
$B^{10,3}=(3,3,1,0,0,\ \dots\ )$}
\end{center}
\vspace{-0,6cm}
\end{figure}

Under the name of parking function, hashing
with linear probing has been and is still
studied by combinatorists
 \cite{FOA,FRA,RIO,SHU,STAN,STAN2}.
Section 4 of    \cite{HASH} contains nice developments
 on the connections
between parking functions and many other combinatorial objects.
In this paper,  we use mainly
  a  -  maybe less exploited  -  connection
between parking functions  and \textit{empirical processes}
of mathematical statistics (see also  \cite{CM,PITSTAN}) .

Let $B_k^{m,\ell}$ denote the size of the $k^{th}$
largest block of consecutive
 cars, and let $B^{m,\ell}= (B_k^{m,\ell})_{k\geq 1}$ be
 the decreasing sequence of
sizes of blocks, ended by an infinite sequence
of $0$'s.
Pittel  \cite{PITT} proves that when $\ell/m$ goes
 to some positive limit $\beta$,
 $B_1^{m,\ell}$ satisfies
\[B_1^{m,\ell}=\frac{2\log m-3
\log\log m+\Xi_m}{2(\beta-1-\log\beta)},\]
 where $\Xi_m$ converges weakly to an
 extreme-value distribution.
This paper is concerned with
  what we would call
the "emergence of a giant block", by reference to
the emergence of a giant component \cite{AGR,ALON,BOLL,FKP,JKLP}.
We have:
\begin{theo}
\label{phase}
For $m$ and $n$ going jointly to $+\infty$
\begin{itemize}
\item[($i$)] if $\sqrt m= o(\ell)$, $B_1^{m,\ell}/m \build
 {\tend}{}{P} 0$;
\item[($ii$)] if $\ell= o(\sqrt m)$, $B_1^{m,\ell}/m \build
 {\tend}{}{P} 1$.
\end{itemize}
\end{theo}
Thus a phase transition occurs for $\ell=\Theta(\sqrt m)$.
The main result of this paper is the description of this phase transition
with the help of Brownian motion theory, following \cite{AGR}.
More precisely, as in  \cite{AGR}, the asymptotic
 behaviour of blocks' sizes is described by widths of excursions of
 stochastic processes related to the Brownian motion. It turns
out, by nature of the problem, and also owing to previous
 works of Aldous \& Pitman \cite{ADD},
that the description  given here (specially by Theorem \ref{lambeaux})
is  more precise than in    \cite{AGR}.

\subsection{Phase transition and Brownian motion}

Recall some notations and definitions from Brownian motion theory.
An \textit{excursion} (from $0$) of the  function
$f$ is the restriction of  $f$
to an interval $[a,b]$  such
that \[f(a)=f(b)=0\textrm{ and }|f(x)|>0\hspace{0,2cm}
 \forall x\in]a,b[;\]
$b-a$ is the \textit{width} or \textit{length} of the excursion, $a$ is the
\textit{starting point} (or the beginning) of the excursion, $b$  the end of
 the excursion. Let us adopt the  notation of \cite[Lecture 4]{YOR} for the
Brownian scaling of a function $f$ over some interval $[a,b]$:
\[f^{[a,b]}=\left(\frac{1}{\sqrt{b-a}}f(a+t(b-a)),
\hspace{0,2cm}0\leq t\leq 1\right).\]
If $f$ is the standard linear Brownian motion, and $g$ (resp. $d$)
is the last zero of $f$ before $1$ (resp. the first zero of $f$ after $1$), then
 $e=\left|f^{[g,d]}\right|$ is called the \textit{normalized Brownian excursion}.
When it is convenient, we regard
the normalized  Brownian excursion
$e(t)$ as defined on the whole real line, being periodic
with period $1$. We define, for $\lambda\ge0$, the operator
$\Psi_{\lambda}$ on the set of bounded functions on the
line by
\bel{psi}
\nonumber\Psi_{\lambda}f(t)&=&f(t)-\lambda t-\inf_{-\infty < s\leq t}(f(s)-\lambda s)\\
 &= &\sup_{s\le t} \left(f(t)-f(s)-\lambda(t-s)\right).
\eel
\begin{figure}[ht]
\begin{center}
\includegraphics[width=7.2cm]{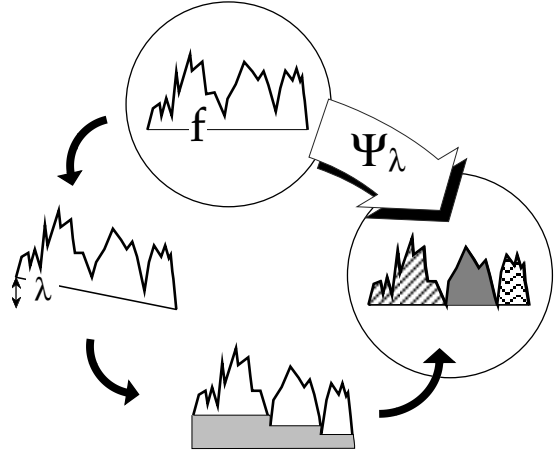}
\end{center}
\caption{$\Psi_{\lambda}f$ and its excursions}
\end{figure}
If $f$ has period 1, then so has  $\Psi_{\lambda}f$.
Evidently, $\Psi_{\lambda} f$ is nonnegative, and we have
\be
\Psi_{\lambda}e(x)&=&e(x)-\lambda x
-\inf_{0\leq y\leq x}(e(y)-\lambda y),\\
\Psi_{\lambda}e(0)&=&\Psi_{\lambda}e(1)=0.
\ee
Let $B(\lambda)=(B_k(\lambda))_{k\geq 1}$
be the sequence of widths
 of excursions of $\Psi_{\lambda}e$,
 sorted in decreasing order.
The  sequence $B(\lambda)$  is a random element of the simplex
\[\{x_1\ge x_2 \ge ... \ge x_n\ge\ ...\,\ge 0,\, \sum_{i\ge 1}x_i=1\}.\]
We have:
\begin{theo}
\label{sorted}
If $\lim \frac{\ell}{\sqrt m}=\lambda\ge 0$,
\[\frac{B^{m,\ell}}{m}\build
 {\tend}{}{law} B(\lambda).\]
\end{theo}

For instance, to complete Theorem
 \ref{phase}, note that
\[B_1^{m,\ell}/m \build{\tend}{}{law} B_1(\lambda).\]
Before we discuss the law of  $B(\lambda)$, in the next Subsection, let us
pursue the description of the asymptotics of the phase transition
 for parking blocks: up to now,
 we only considered the parking process frozen at a given time $n=m-\ell$,
that is, just after the arrival of car $c_n$. The next Theorem describes
 the evolution of blocks' sizes, as cars arrive, during the phase transition:
asymptotically, the joint law of sequences of blocks' sizes
$B^{m,\lceil\lambda_i\sqrt m\rceil}$
after  successive arrivals of cars $c_{n(i)}$,
$n(i)=m-\lceil\lambda_i\sqrt m\rceil$,
 $\lambda_1>\lambda_2>\cdots>\lambda_k$, once these sequences are
normalized,
converges to the joint law of $\left(B(\lambda_i)\right)_{i=1,2, ... ,k}$.
More formally, set
\be
B^{(m)}(\lambda)&=&
\frac{B^{m,\lceil\lambda\sqrt m\rceil}}{m}.
\ee
We have:
\begin{theo}
\label{lambeaux}
The finite-dimensional distributions of $\left(B^{(m)}(\lambda)
\right)_{\lambda\ge0}$
 converge weakly to the finite-dimensional
distributions of  $\left(B(\lambda)
\right)_{\lambda\ge0}$.
\end{theo}
Though, in the random
 graph model, the asymptotic distribution of sizes of clusters
(connected components)  has a description similar to that given at Theorem \ref{sorted},
the analog of  Theorem \ref{lambeaux} is false, as observed by
Aldous \cite{AGR}: in  coalescence models
 based on excursions
of stochastic processes,   clusters (excursions) can only merge
 with their neighbors,
while this is not true for connected components
of  the random graph.
In Section \ref{uniform}, at the price of
heavier notations, we give the
analog of Theorem \ref{lambeaux} for  the
asymptotic behaviour of  sizes \textit{and positions} of blocks.

\subsection{Size-biased  permutations}
\label{sbp}
As a consequence of  \cite[Theorem 4]{PAV}, we have
\begin{theo}
\label{Pavlov}
The distribution function
$\Pr(B_1(\lambda)\le x)$ has the following expression:
\[1+\lambda^3e^{\lambda^2/2}
\sum_{k\ge 1}\frac{(-1)^k}{k!}\int_{D(\lambda,x,k)}
\frac{\lambda^{2k}\exp\{-\lambda^4/2(\lambda^2-
x_1-\dots-x_k)\}dx_1\dots dx_k}{(2\pi)^{k/2}(x_1\dots x_k(
\lambda^2-x_1-\dots-x_k))^{3/2}},\]
in which
\[D(\lambda,x,k)=\left\{(x_i)_{1\le i\le k}\,:\ x_i\ge \lambda^2x,
1\le i\le k,\textrm{ and }\sum x_i\le \lambda^2\right\}.\]
\end{theo}
Theorem 4  of  \cite{PAV}  gives the limit law of the  largest tree in a random forest:
it turns out that  forests and parking schemes are in  one-to-one
correspondence  (see Subsection \ref{prPavlov}).
Flajolet \& Salvy \cite{FLAS} have a direct approach, to the
computation of the density of  $B_1(\lambda)$, by methods based on Cauchy
 coefficient integrals to which the saddle point method is applied: the density
 they obtain is a variant of the Dickman function \cite[Ch. III, Sec. 5.3]{TEN}.

In view of Theorem \ref{Pavlov},  the joint law of
 $(B_1(\lambda),B_2(\lambda),\,\dots\,,B_k(\lambda))$ seems out of reach,
but we are more lucky with the joint law of the first terms of a sequence
 $R(\lambda)$  obtained by permutation of the terms of  $B(\lambda)$.
Roughly speaking, in the \textit{size-biased
 permutation} $R(\lambda)$ of a random probability distribution such as
$B(\lambda)$, the largest terms of the sequence $B(\lambda)$  appear with a high
probability  at the beginning of the sequence $R(\lambda)$:  we have
\begin{equation}
\label{coupon}
\Pr\left(R_1(\lambda)=B_{k}(\lambda)\ \vert\ B(\lambda)\right)
=B_{k}(\lambda),
\end{equation}
the $k^{th}$ term of $R(\lambda)$  being also drawn randomly
  with a probability proportional to its size,
but among the terms  that  did not appear before.
A more formal  definition of size-biased
 permutations, by construction through a rejection method,   is given in \cite{PIT2}:
 consider a sequence of independent, positive, integer-valued random variables
$(I_k)_{k\geq 1}$, distributed according to  $B(\lambda)$:
\[\Pr\left(I_k=j\ \vert\ B(\lambda)\right)=B_j(\lambda).\]
With probability $1$ the terms of $B(\lambda)$ are positive, as
 $\Psi_{\lambda}e$ has infinitely many excursions, so each positive
 integer appears at least once in the sequence $(I_k)_{k\geq 1}$. Erase
 each repetition after the first  occurence  of a given integer in the sequence: there
remains a random permutation $(\sigma(k))_{k\geq 1}$ of the positive integers.  Set:
\begin{equation}
\label{randperm}
R_k(\lambda)=B_{\sigma(k)}(\lambda).
\end{equation}
We have:

\begin{theo}
\label{mixed}
The law of the size-biased
 permutation  $R(\lambda)$ of
$B(\lambda)$
satisfies
\[\left(R_1(\lambda)+R_2(\lambda)+~~...~~+
R_k(\lambda)\right)_{k\geq 1}
\build
 {=}{}{law} \left(\frac{N_1^2+N_2^2+~~...~~+N_k^2}
{\lambda^2+N_1^2+N_2^2+~~...~~+N_k^2}\right)_{k\geq 1},\]
in which the $N_k$ are standard Gaussian and independent.
\end{theo}

Actually,
Theorem \ref{mixed}
 gives an implicit description of the law
of $B(\lambda)$, for instance it proves that almost surely each
$R_k(\lambda)$ is positive, and thus a.s. $0<B_k(\lambda)<1$.
Size-biased permutations of random discrete probabilities
have been studied, among others, by  Aldous \cite{FLOUR}
 and Pitman \cite{PIT1,PIT2}.
The most celebrated example is the size-biased permutation
of the sequence of limit sizes of
cycles of  a random permutation.
While the limit distribution of the sizes
of the largest, second largest ...
cycle have a complicated expression
 \cite{DICK,SHLL},
the successive terms $R_1,R_2,~~...~~$ of their size-biased
 permutation satisfies
\[(R_1+R_2+~~...~~+
R_k)_{k\geq 1}
\build
 {=}{}{law} \left(1-U_1U_2~~...~~U_k\right)_{k\geq 1},\]
in which the $U_k$ are uniform on $[0,1]$ and independent.
Actually,  it  is common that the distribution of the size-biased
 permutation  of a sequence has a simpler distribution than the original sequence,
when the sequence  is
 related  to a Poisson point process, a famous example being the Poisson-Dirichlet
 distribution  \cite{ABT1,ABT2,PER,PPY,PITYOR}.
The distribution of  $R(\lambda)$, as described in Theorem \ref{mixed},
already appeared as the  law of the
$\Delta$\textit{-valued fragmentation
 process}
derived from the continuum random
 tree, introduced by Aldous \& Pitman
in their study of the standard additive
coalescent \cite[Corollary 5]{ADD}:
this is commented in the next Subsection.

As in the case of sizes of cycles, the  unnatural size-biased permutation
of   $B(\lambda)$  is the limit of a natural permutation of
$B^{(m)}(\lambda)$: define  $R^{m,\ell}=
\left(R_k^{m,\ell}\right)_{k\geq 1}$
 as the sequence of sizes of blocks
 when the blocks are sorted by increasing
 date of birth
(in increasing order of first arrival
 of a car). If $\ell\geq m$, or if there are less than $k$ blocks,
 set $R_k^{m,\ell}=0$.
For instance, on
Figure 1, $B^{m,\ell}=(3,3,1,0,\ \dots\ )$ and
$R^{m,\ell}=(3,1,3,0,\ \dots\ )$. Concerning $R^{m,\ell}$, we have an analog of Theorem
\ref{sorted} :
\begin{theo}
\label{size-biased}
If $\lim m^{-1/2}\ell=\lambda\ge 0$,
\[\frac{R^{m,\ell}}{m}\build
 {\tend}{}{law} R(\lambda).\]
\end{theo}
Set
\be
R^{(m)}(\lambda)&=&
\frac{R^{m,\lceil\lambda\sqrt m\rceil}}{m}.
\ee

\begin{figure}[ht]
\begin{center}
\includegraphics[width=7.2cm]{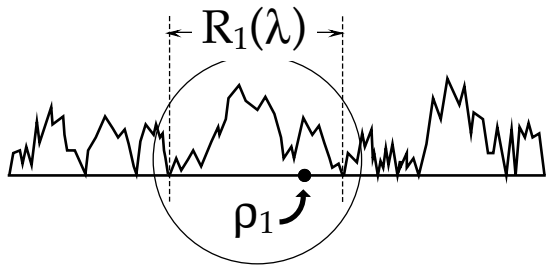}
\end{center}
\caption{Random generation of $R_1(\lambda)$}
\label{r1delambda}
\end{figure}

For an analog of
Theorem \ref{lambeaux} to hold true, giving the convergence
 of finite dimensional distributions of $R^{(m)}$, we should define
 $\left(R(\lambda)\right)_{\lambda\ge0}$
 as a process. This is not straightforward, as there are
many possible definitions  of the size-biased permutation
$\sigma_{\lambda}$ that  throws $B(\lambda)$ on $R(\lambda)$ :
$\sigma_{\lambda}$  has to be defined as a process too. For sake of brevity, we
 shall only state a result for the first component of  $R(\lambda)$.
Consider a random number $\rho_1$, uniform on $[0,1]$ and independent of $e$,
and define $R_1(\lambda)$ as the width of the excursion of $\Psi_{\lambda}e$
that contains $\rho_1$  (see Figure \ref{r1delambda}). Then $R_1(\lambda)$ 
is defined, and satisfies (\ref{coupon}),
simultaneously for each value of $\lambda$. We have:
\begin{theo}
\label{fidis2}
The finite-dimensional distributions of $\left(R_1^{(m)}(\lambda)
\right)_{\lambda\ge0}$
 converge weakly to the finite-dimensional
distributions of  $\left(R_1(\lambda)
\right)_{\lambda\ge0}$.
\end{theo}
It turns out that the limit process has a rather simple description: set
\[\Sigma(\lambda)=-1+\frac{1}{R_1(\lambda)}.\]
We have
\begin{theo}
\label{Levy}
 $\Sigma=(\Sigma(\lambda))_{\lambda\geq 0}$
is the stable subordinator
with exponent $1/2$.
\end{theo}
It is well known that the family of first hitting
times of levels $\lambda$ by the Brownian motion is also  the stable subordinator
with exponent $1/2$. The following description of the stable subordinator
with exponent $1/2$, by its finite dimensional distributions, will be useful
for the proof: for any $k$ and any k-tuple of positive numbers $(\lambda_i)
_{1\leq i \leq k}$,
\begin{equation}
\label{fidis}
\left(\Sigma(\lambda_1+\lambda_2+ ... +\lambda_i)\right)_{1\leq i \leq k}
\build{=}{}{law}\left(\frac{\lambda_1^2}{N_1^2}+
\frac{\lambda_2^2}{N_2^2}+
 ... +\frac{\lambda_i^2}{N_i^2}\right)_{1\leq i \leq k},
\end{equation}
in which the $N_k$ are standard Gaussian and independent.

The well known  fact that $\Sigma(\lambda)$ is a pure jump process makes
sense in the parking
 scheme context, since the block of car $c_1$ is known to increase by $O(m)$
while only $O(\sqrt m)$ cars arrived: it can only be explained by coalescence
with other blocks of size $O(m)$, that is, by \textit{instantaneous jumps}.
Incidentally, let $L(\lambda)$ denote the length
of the excursion of $\Psi_{\lambda}e$ beginning at $0$, and set
\[\tilde{\Sigma}(\lambda)=-1+\frac{1}{L(\lambda)}.\]
  Bertoin \cite{BERT}  nicely proves that $(\Sigma,R_1)$ and $(\tilde{\Sigma},L)$
 have the same law.
For the moment, we do not see any combinatorial explanation of this
 identity between $R_1$ and $L$.

\subsection{Coalescence}
\label{coal}

We give here a brief account of coalescence,  which is masterfully surveyed
in  \cite{ALD2,ALD3}. We essentially quote the two previously cited references.
Models of  coalescence (aggregation, coagulation, gelation ...) have been studied
in many  scientific disciplines, essentially physical chemistry, but also astronomy,
bubble swarms, mathematical genetics, and recently random graph theory
\cite[Section 1.4]{ALD3}. In a basic model, clusters with different masses
move through space, and when two clusters (say, with masses $x$ and $y$)
are sufficiently close, there is some chance that they merge into a single cluster
of mass $x+y$ \cite[Section 1.1]{ALD3}. The probability
that they merge is quantified, in some sense, by a \textit{rate kernel}
$K(x,y)$.  As far as parking is concerned, the growth of clusters (parking blocks) is due
partly  to cars' arrivals, partly to aggregation with other blocks,
but  we saw that during the phase transition  the coalescence factor is preponderant.

A complete model for coalescence, detailing mass, position, and velocity of each cluster,
is too complicated for analysis, so recent works focused on
 the evolution of masses of clusters through time: the \textit{general stochastic coalescent}
 \cite{EPIT}
is the continuous-time Markov process whose state space is the infinite-dimensional
simplex \[\Delta=\left\{(x_i)_{i\ge 1}\ :\ x_i \ge 0,\ \sum x_i=1\right\},\]
(the $x_i$'s are the sizes of clusters) and that evolves according to the rule
\[\textrm{each pair }(x_i,x_j)\textrm{ of clusters  merges at rate }K(x_i,x_j).\]
It means  that, if at time $t$ the state of the system is $(x_i)_{i\ge 1}$, the next pair $(I,J)$
of clusters that will merge and the time $t+T$ when they merge
 are jointly distributed as follows:
assume we are given a set of independent random variables $(T_{i,j})_{1\le i<j}$
with distribution described by
\[\Pr\left(T_{i,j}>t\right)=\exp\left(-K(x_i,x_j)t\right),\]
and set
\[\inf_{1\le i<j}T_{i,j}=T_{I,J}=T.\]

 It turns out that the way  connected components merge in the random
graph process is somehow related to the \textit{multiplicative coalescent}
($K(x,y)=xy$) \cite{AGR}. One rather expects the parking to be related to
 the \textit{additive coalescent} ($K(x,y)=x+y$): given that a parking scheme
 with  $m$ places, $n$ cars and  $\ell=m-n$ empty places has two blocks with
 size $x$ and $y$, the probability that these two blocks merge at the next arrival is
\begin{equation}
\label{pt}
\frac{x+y+2}{(\ell-1)m},
\end{equation}
as the number of empty places after block $x$ but before block $y$ is random uniform
on $1,2,\dots, \ell-1$, and, given that this number is $1$ (resp. $\ell-1$, $\notin
 \{1,\ell-1\}$) the conditional probability that the two blocks merge at the next arrival is
$\frac{x+1}{m}$ \B(resp. $\frac{y+1}{m}$, $0$\B).
Aldous \& Pitman \cite{ADD} give a construction of the additive coalescent
 through a fragmentation process $Y=\left(Y(\lambda)\right)
_{\lambda\ge 0}$: the $\Delta$-valued random variable $Y(\lambda)$ is the  ranked
sequence of masses  of tree components of continuum forests obtained by cutting
 the "edges" of the \textit{Brownian continuum random tree}
by a Poisson process of cuts with rate $\lambda$ by unit length. As more or less
expected, according to Theorem  \ref{mixed}   and to \cite[Corollary 5]{ADD},
the distributions of $Y(\lambda)$  and  $B(\lambda)$ are the same.

Furthermore, let $\rho_1^{\star}$ be a leaf, of the
Brownian continuum random tree, picked uniformly at random according to the
mass measure, and let $Y_1^{\star}(\lambda)$ denote the mass of the tree component
of the random forest
that contains  $\rho_1^{\star}$  when the cutting intensity is $\lambda$. Then, according
to Theorem \ref{Levy} and to \cite[Theorem 6]{ADD}, the distributions of the
stochastic processes $Y_1^{\star}$ and  $R_1$ are the same. These facts
suggest that
\begin{theo}
\label{id}
The \textit{processes} $B$ and $Y$ have the same distribution.
\end{theo}
Theorem \ref{id}  is actually the main result of a recent paper by Bertoin
\cite{BERT}. In Section \ref{addcoal},
 we give an alternative proof of Theorem \ref{id},  that relies on
Theorem \ref{decomp},
 a path decomposition result for $\Psi_{\lambda}e$.

 Note that Theorem \ref{Levy} is not a mere consequence  of Theorem \ref{id}
and \cite[Theorem 6]{ADD}, as the very similar selection mechanisms leading to
$Y_1^{\star}$ (resp.  $R_1$) depend not only on the stochastic
processes $Y$ (resp.  $B$), but on underlying richer structures,
a family of Poisson point processes of cuts of a Brownian continuum
 random tree on one hand,  and the  family of stochastic processes
$\Psi_{\lambda}e$ on the other hand. Even if one of the  constructions of the
 Brownian continuum random tree uses the normalized Brownian excursion
 \cite[Corollary 22]{CRT3}, we do not know for the moment any extension of
Theorem \ref{id} to these richer structures, that
 would yield a direct proof of the identity
between the distributions of the  stochastic processes
 $Y_1^{\star}$ and  $R_1$.
However, in the concluding remarks, we give a rather
 convincing combinatorial explanation
of the connection between the two richer structures.

\subsection{Decomposition of sample paths of $\Psi_{\lambda}e$}
\label{paths}

In previous subsections, objects from Brownian motion theory
allowed to describe phase transition for parking schemes. In this subsection, we translate
 the parking schemes combinatorial identity:
 \[m^{n}=\sum_{k=1}^{n}C_{n-1}^{k-1}m(k+1)^{k-1}
(m-k-1)^{n-k-1}(m-n-1)\]
to obtain Theorem \ref{decomp},
 a property of decomposition of sample paths of $\Psi_{\lambda}e$  used
 in Section \ref{addcoal}
  to give  simple proofs  of   Theorems  \ref{mixed},  \ref{Levy} and \ref{id}.
\begin{figure}[ht]
\begin{center}
\includegraphics[width=8cm]{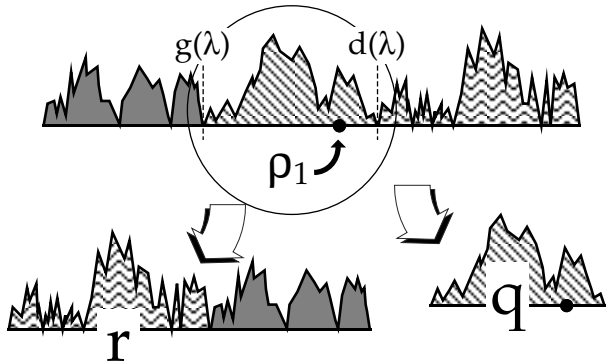}
\end{center}
\label{thegap}
\caption{Decomposition  of $\Psi_{\lambda}e$.}
\end{figure}
Let $\rho_1$ be  a random variable uniformly distributed
 on $[0,1]$ and independent of $e$. Almost surely, $\Psi_{\lambda}e(\rho_1)$ is positive.
Let $g(\lambda)$ (resp. $d(\lambda)$) denote the last
 zero of $\Psi_{\lambda}e$ in the interval $[0,\rho_1)$
(resp. the first zero in the interval $(\rho_1,1]$),
so that $R_1(\lambda)=d(\lambda)-g(\lambda)$.
To avoid the extensive use of notation
$\{x\}$ for the fractional part of the real number $x$, we shall extend
 $\Psi_{\lambda}e$, as well as other functions defined on $[0,1]$, such as $q$ or $r$
defined below, to  periodic
 functions  on the line. We set
\begin{eqnarray*}
q&=&(\Psi_{\lambda}e)^{[g(\lambda),d(\lambda)]},\\
r&=&(\Psi_{\lambda}e)^{[d(\lambda),g(\lambda)+1]}.
\end{eqnarray*}
Let $\tau_x$ denote the shift operator for functions  on the line,
defined by
\[\left(\tau_xf\right)(y)=f(x+y).\]
\begin{theo}
\label{decomp}
 We have:
\begin{itemize}
\item[(i)] $R_1(\lambda)$ has the same distribution
as $\frac{N^2}{\lambda^2+N^2}$,
 in which $N$ is standard Gaussian ;
\item[(ii)] $q$
is a normalized Brownian excursion, independent of  $R_1(\lambda)$ ;
\item[(iii)] Let $w$ be uniformly distributed on $[0,1]$ and independent of $e$.
 Given $(q,\rho_1)$ and $R_1(\lambda)=x$, $\tau_wr$ has the same distribution
as $\tau_w\Psi_{\frac{\lambda}{\sqrt{1-x}}}e$.
\end{itemize}
 \end{theo}

Actually, not only the conditional distribution of  $\tau_wr$, but also
the conditional distribution of  $r$ has a simple description in terms of the
  Brownian motion, and also as a nonuniform random shift of
  $\Psi_{\frac{\lambda}{\sqrt{1-x}}}e$ \cite{CJ}.  However, the weaker form $(iii)$
fills our needs for the proofs of
  Theorems  \ref{mixed},  \ref{Levy} and  \ref{id}.

\begin{figure}[ht]
\hspace{71pt}
\includegraphics[width=9.05cm]{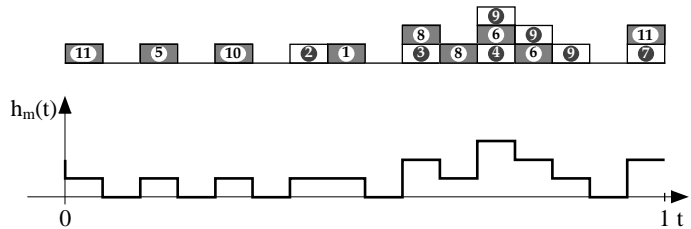}
\caption{A parking scheme and its profile ($m=16, n=11$)}
\label{profileps}
\end{figure}

The paper is organized as follows.
Section \ref{given} analyses the block containing
a given car or a given site, leading to the proof of Theorem
 \ref{phase}. At  Section \ref{profil}, we give the proof of the main result,
 Theorem \ref{sorted},
 with the help of a close
 coupling between \textit{empirical processes} of mathematical
statistics and the \textit{profile} obtained  by assuming that each car
lays a $1/\sqrt m$-thick layer of sediment on the way between its first try
and its final place
(see Figure \ref{profileps}).
We extend these arguments at Section \ref{uniform}
to obtain the asymptotic of the joint law,
at different times, of widths and \textit{positions}
of blocks.
Distributional results, Theorems \ref{Pavlov}  and \ref{size-biased}, are proven
at Section \ref{Markov} by combinatorial arguments.
We prove  Theorems \ref{fidis2} and \ref{decomp} at Section \ref{fractal},
with the help of Theorem \ref{monticule}, about weak convergence of profiles.
Finally, in  Section \ref{addcoal},   Theorems \ref{mixed},
 \ref{Levy} and  \ref{id} are shown to be consequences of Theorem \ref{decomp}.
 Section \ref{concl} concludes the paper with an attempt of combinatorial explanation
 for the connections between our paper and \cite{ADD}.

\section{On the block containing
 a given car, or a given site}
\label{given}

In this Section, we prove Theorem
 \ref{phase}, with the help  of a weaker form of Theorem
 \ref{size-biased}, concerning the
size $R^{m,\ell}_1$ of the block containing
car $c_1$: we have
\begin{theo}
\label{undemi}
If $m^{-1/2}\ell\tend\lambda>0$,
\[\frac{R^{m,\ell}_1}{m}\hspace{0,3cm}
\build{\tend}{}{law}\hspace{0,3cm}
 \frac{N^2}{\lambda^2+N^2},\]
in which $N$ is standard Gaussian.
\end{theo}

\pr  Let $f(\lambda,x)$ denote the
 density of $\frac{N^2}{\lambda^2+N^2}$,  and let
$\varphi(m,n,k)$ denote the probability  that,
when parking $n=m-\ell$ cars on $m$ places,
the block containing car $c_1$ has $k$ elements. We have
\begin{eqnarray}
\nonumber f(\lambda,x)&=&
\frac{\lambda}{\sqrt{2\pi}}\ x^{-1/2}(1-x)^{-3/2}
\exp\left(-\frac{\lambda^2x}{2(1-x)}\right)
1_{]0,1[}(x),\\
\label{combi}\varphi(m,n,k)&=&C^{k-1}_{n-1}
\frac{(k+1)^{k-1}}{m^n}\ m(m-k-1)^{n-k-1}
(\ell-1).
\end{eqnarray}
From the  change of variable
 $x=y/(\lambda^2+y)$, leading to
\be
\int_0^1f(\lambda,x) dx=
\int_0^{+\infty}\frac{e^{-y/2}dy}{\sqrt{2\pi y}},
\ee
we deduce that $f(\lambda,x)$ is a density of probability
and that if some random variable $X$
has the density $f(\lambda,x)$,
 then $\lambda^2X/(1-X)$ has a
 $\gamma_{1/2,1/2}$ law, the  law of the square of
 a standard Gaussian random
 variable.
That is, $N^2/(\lambda^2+N^2)$
 has density $f(\lambda,x)$.

To explain (\ref{combi}), first we remark that the number of parking schemes
for $n$ cars on $m$ places is  $m^n$. If we specify that the last place has to be empty,
we get what is called a \textit{confined} parking scheme: there are $(m-n)m^{n-1}$
confined parking schemes \cite{HASH,KNU2}, as each orbit drawn by the group of rotations
has $m$ elements, among which $m-n$ are confined. A block with $k$ cars can be seen
as a confined parking scheme of $k$ cars on $k+1$ places, so there are
  $(k+1)^{k-1}$ ways to build such a block.
Turning to  (\ref{combi}), one has to choose the set of $k-1$ cars that belong
 to the same block as $c_1$, giving the factor $C_{n-1}^{k-1}$,
the place where this block begins, giving the factor $m$, the way
 these $k$ cars are allocated on these $k$ places, giving the
factor $(k+1)^{k-1}$, and finally one has to park
the $n-k$ remaining cars on the $m-k-2$ remaining places,
leaving one empty place at the beginning and at the end of the block
containing car $c_1$. This can be done in $(m-k-1)^{n-k-1}(\ell-1)$ ways,
 the number of \textit{confined} parking schemes of $n-k$ cars on
$m-k-1$ places. Note
that these computations would hold for any given car instead of $c_1$.

For $0<a<b<1$  and $m^{-1/2}\ell\tend\lambda$,
\[\lim_m\ \Pr(am\leq R^{m,\ell}_1\leq bm)
=\int_a^bf(\lambda,x) dx,
\]
is a straightforward consequence of
\begin{lem} For any $0<\varepsilon<1/2$
there exists a constant $C(\varepsilon)$
such that, whenever, simultaneously, $\varepsilon\leq\frac{k}{m}
\leq 1-\varepsilon$ and $\varepsilon\leq \frac{\ell}{\sqrt m}
\leq \frac{1}{\varepsilon}$, we have:
\[\left|\varphi(m,n,k)-\frac{1}{m}
f\left(\frac{\ell}{\sqrt m},\frac{k}{m}\right)\right|\leq
C(\varepsilon)m^{-3/2}.\]
\label{local}
\end{lem}
Lemma  \ref{local} is proven at
 the end of this Section. $\hspace{1cm}\diamondsuit$

\textit{Proof of Theorem} \ref{phase}$(ii)$.
We assume  $\ell=o(\sqrt m)$.
Provided that $\ell\leq\lambda\sqrt m$,
\[B_1^{m,\ell}\geq mR^{(m)}_1(\lambda).\]
Thus, for any $\lambda>0$
 and for $m$ large enough:
\[\Pr(B_1^{m,\ell}<mx)\leq
\Pr\left(R^{(m)}_1(\lambda)<x\right).\]
Due to Theorem \ref{undemi},
we obtain that for any $\lambda>0$
\[\limsup_m\Pr(B_1^{m,\ell}<mx)\leq
\Pr\left(\frac{N^2}{\lambda^2+N^2}<x\right).\]
Clearly, for $x<1$,
\[\inf_{\lambda>0}\Pr\left(\frac{N^2}{\lambda^2+N^2}<x\right)=0.\]

\textit{Proof of Theorem} \ref{phase}$(i)$.
Let $L^{(m)}(\lambda)$
 be the length, normalized by $m$, of the block
of cars containing \textit{place} $1$
when car
$c_{\lfloor m-\lambda\sqrt m\rfloor}$
has parked. We have, for $k>0$,
\[\Pr\left(L^{(m)}(\lambda)=\frac{k}{m}\right)=
\frac{\lfloor m-\lambda\sqrt m\rfloor}{m}
\Pr\left(R^{(m)}_1(\lambda)=\frac{k}{m}\right),\]
and place $1$ is empty with probability:
 \[\Pr\left(L^{(m)}(\lambda)=0\right)=
\frac{\lceil \lambda\sqrt m\rceil}{m}.\]
We also have
\[\Pr\left(L^{(m)}(\lambda)=\frac{k}{m}\ \left\vert\ B^{(m)}_1
(\lambda)=\frac{k}{m}\right.\right)\geq\frac{k}{m},\]
and thus
\begin{eqnarray*}
E\left[L^{(m)}(\lambda)\right]&\geq&\sum_k
\frac{k}{m}\Pr\left(L^{(m)}(\lambda)=B^{(m)}_1(\lambda)=\frac{k}{m}\right)
\\
&\geq&E\left[\left(B^{(m)}_1(\lambda)\right)^2\right].
\end{eqnarray*}
Owing to $\sqrt m=o(\ell)$, we obtain that for any
$\lambda>0$,
\[B_1^{m,\ell}(\omega)
\leq mB^{(m)}_1(\lambda,\omega),\]
when $m$ is large enough, not depending on $\omega$, so that:
\begin{eqnarray*}
\limsup_mE\left[\left(\frac{B^{m,\ell}_1}{m}\right)^2\right]&\leq&
\inf_{\lambda>0}\hspace{0,2cm}
\lim_m\ E\left[L^{(m)}(\lambda)\right]\\
&=&\inf_{\lambda>0}\hspace{0,2cm}
E\left[\frac{N^2}{\lambda^2+N^2}\right],
\end{eqnarray*}
 yielding  $(i)$.$\hspace{1cm}\diamondsuit$

\textit{Proof of Lemma \ref{local}}. Setting, for brevity,
$x=k/m$ and $\lambda=m^{-1/2}\ell$, we can write
\[\varphi(m,m-\lambda\sqrt m,k)=\Phi_{m,1}\Phi_{m,2}\Phi_{m,3},\]
in which:
\be
\Phi_{m,1}&=&C^{xm-1}_{m-\lambda \sqrt m-1}\\
\Phi_{m,2}&=&(xm+1)^{xm-1}m^{-m+\lambda \sqrt m-1}\\
\Phi_{m,3}&=&(m-xm-1)^{m-\lambda \sqrt m-xm-1}
(\lambda \sqrt m-1).
\ee

We obtain
\be
\Phi_{m,1}&=&\frac{x^{-xm+1/2}e^{\kappa(m)-\beta(m)}
\left(1+O(1/k)+O\left(\sqrt m/(m-k)\right)\right)}
{\sqrt{2\pi m}\hspace{0,5cm}
(1-x)^{(1-x)m-\lambda\sqrt m+1/2}}\\
\kappa(m)&=&1+(m-\lambda \sqrt m-1/2)\log(1-\lambda/\sqrt m-1/m)\\
&=&-\lambda \sqrt m+\lambda^2/2+O(m^{-1/2})\\
\beta(m)&=&(m-\lambda \sqrt m-xm+1/2)\log\left(1-
\frac{\lambda}{(1-x)\sqrt m}\right)\\
&=&-\lambda \sqrt m
+\frac{\lambda^2}{2(1-x)}
+O\left((m-k)^{-1/2}\right).\\
\Phi_{m,2}&=&x^{xm-1}m^{xm-m+\lambda \sqrt m}
(1+O(1/k))e\\
\Phi_{m,3}&=&\frac{\lambda m^{(1-x)m-\lambda \sqrt m}
(1-x)^{(1-x)m-\lambda \sqrt m}}{e(1-x)\sqrt m}
\left(1+O\left(\sqrt m\right)+O\left(\frac{\sqrt m}{m-k}\right)\right)
\ee
and finally:
\be
\varphi(m,m-\lambda\sqrt m,k)&=&
\frac{\lambda(1-x)^{-3/2}}{m\sqrt{2\pi x}}
\exp\left(-\frac{\lambda^2x}{2(1-x)}\right)(1+\eta(m,k))\\
|\eta(m,k)|&\leq&\frac{K_1}{k}+\frac{K_2\sqrt m}{m-k}
+\frac{K_3}{\sqrt m}
+\frac{K_4}{\sqrt{m-k}}.\hspace{1cm}\diamondsuit
\ee

\section{Profiles of  parking schemes}
\label{profil}
Let $H_k$ denote the number of cars that
tried to park on place $k$,
\textit{successfully or not}, and let  $h_m$ denote the
\textit{profile} of the parking scheme, defined by:
\be
h_m(t)&=&\frac{H_{\lfloor mt\rfloor}}{\sqrt m}.
\ee
As $H_k=0$ if and only if place $k$ is empty,
 the width of an excursion
of $h_m$ turns out to be the length
 of some block of cars, normalized by
$1/m$.
Set:
\be
h_{\lambda}(t)&=&\Psi_{\lambda}e(\{-v+t\}),
\ee
in which $v$ denotes a uniform random variable
independent of $e$.

In this Section we give the proof of Theorem \ref{sorted},  that has roughly speaking
three steps: as a first result, we
 establish in Subsection \ref{connect}
a close coupling
 between $H_k$ and the empirical
 processes of mathematical statistics.
In  Subsection \ref{exc},
using Theorems of Donsker and Vervaat,
 we  prove the following  Theorem, which is the key to this paper.
\begin{theo}
\label{monticule}
If $\lim_m\frac{\ell}{\sqrt m}=\lambda$,
\[h_m\build{\tend}{}{weakly}h_{\lambda}.\]
\end{theo}
Theorem \ref{sorted} states the convergence of  widths of
 excursions of $h_m$ to widths of
 excursions of
 $h_{\lambda}$. Its proof, given in  Subsection \ref{ald97},
  requires some care, as the sequence of widths of excursions
 is not a continuous functional of $h_m$:
 the proof relies on  an extension of the invariance principle
that we learned from   \cite{AGR,APIT}.
Further consequences
 of Theorem \ref{monticule} are
Theorem \ref{decomp} and also some results
 about stochastic processes
developped in \cite{CJ}. Theorem \ref{lambeaux} is the consequence
 of Theorem \ref{2dim},
an extension of  Theorem \ref{monticule}. The case
$\lambda=0$, $\ell=1$ of  Theorem
\ref{monticule} was developped
in \cite[Section 4]{CM} for the study of
the width of labeled trees.

\subsection{Connection between
 parking and empirical processes}
\label{connect}

Propositions \ref{Skoro}
and \ref{record},
 at the end of this subsection, are the key points
for the convergence of blocks' sizes.
Given a sequence $(U_k)_{k\geq 1}$ of
independent uniform random variables,
 we assume that first
try of car  $c_k$ is place $\lceil mU_k\rceil$, all parking schemes
being thus equiprobable.
If $Y_k$ denotes the number of cars whose first
try was  place $k$, then we have:
\begin{equation}
H_{k+1}=Y_{k+1}+(H_k-1)_+,
\label{reccirc}
\end{equation}
since either place $k$ is occupied by car $c_i$
 and, among the $H_k$ cars
 that tried place $k$, only car
 $c_i$ won't visit place $k+1$, so that
$H_{k+1}=Y_{k+1}+(H_k-1)$, or
 place $k$ is empty and $H_{k+1}=Y_{k+1}$.
We understand this equation,
when $k=m$, as $H_1=Y_1+(H_m-1)_+$.
This induction alone does not
give the $H_k$'s, since we do not have any
starting value.
 The gap is filled by Proposition \ref{empty}, that gives the
 connection between hashing (or parking)
 and empirical processes.

Given a sample $(U_1,U_2, ...,U_n)$
of uniform random variables, the
\textit{empirical distribution} $F_n$ and the
\textit{empirical process} $\alpha_n$
(see
\cite{CSO,POLL,SHOW}
for background)
are respectively defined by
\be
F_n(t)&=&\frac{\#\{1\leq i\leq n\hspace{0,1cm}|
\hspace{0,1cm}U_i\leq t\}}{n}\\
&=&t+\frac{\alpha_n(t)}{\sqrt n}.\ee
The process $\alpha_n$ gives a measure of the accuracy
of the approximation of the true distribution function $t$
by the empirical distribution function
 $F_n(t)$, and was, as such,
extensively studied in mathematical statistics.
Let $V$ be defined  by
\be
a&=&\min\{\alpha_n(k/m)\ |\ 1\leq k\leq m\}\\
V&=&\min\{j\ |\ 1\leq j\leq
 m\textrm{ and }\alpha_n(j/m)=a\}.
\ee

\begin{pro}
Place $V$ is empty.
\label{empty}
\end{pro}
\textit{Proof of Proposition \ref{empty}}.
Set:
\be
A_k&=&\sqrt n\ \alpha_n(k/m)\\
&=&\#\{1\leq i\leq n\hspace{0,1cm}|
\hspace{0,1cm}U_i\leq k/m\}-k
\hspace{0,1cm}\frac{n}{m}.
\ee
Since we have:
\[Y_k=\#\{j\ |\ 1\leq j\leq n\textrm{ and }\lceil mU_j\rceil=k\},\]
it follows that:
\begin{equation}
A_{k+1}=A_k+Y_{k+1}-n/m.
\label{rec}
\end{equation}
As  $A_0=A_m$, we can extend $A_k$,
and $Y_k$ as well,
to  periodic sequences, so that (\ref{rec})
holds true for any integer.
Thus
\begin{eqnarray}
\nonumber  Y_{V-k+1}+Y_{V-k+2}+\cdots+Y_{V}
&=&kn/m-A_{V-k}+A_{V}\\
\label{plein}&\le&\lfloor kn/m\rfloor\le k-1,
\end{eqnarray}
the first inequality by definition of $V$.
We remark  that if the set of places
$\{j-k+1,j-k+2,\cdots,j\}$
is full while no more than $k-1$ cars had their first try in it,
then necessarily place $j-k$
is occupied.
Thus if place $V$ is not empty,  by letting $k=1$ in (\ref{plein})
we obtain  that $V-1$ is not empty either,
and by induction on $k$, still using (\ref{plein}), no place is  empty.
$\hspace{1cm}{\diamondsuit}$

For $0\le k\le m$, set:
\be
C_k&=&\sqrt n\ \alpha_n(k/m)-k
\hspace{0,1cm}\frac{\ell}{m}\\
&=&Y_{1}+Y_{2}+\cdots+Y_{k}-k,
\ee
and extend it to any integer, through
$C_{k+m}=C_k-\ell$.
With the convention that $H_k$ is periodic as well,
(\ref{reccirc}) holds true for any integer,
and we can use it  to compute  $H_k$,
starting from  $H_{V}=0$:
\begin{pro}
\label{Skoro}
For any $k \in \{1,2,\ ...\ ,m-1\}$,
\[
H_{V+k}=C_{V+k}-C_{V}+
\max_{1\leq i\leq k}\{(C_{V-1}-
C_{V+i-1})_+\}.
\]
\end{pro}
Since blocks of cars are just blocks of consecutive indices $k$
such that $H_k>0$ (\textit{excursions} of $H_k$), our study relies essentially
on this expression, that connects blocks of cars with empirical processes. A
 similar line of proof is used in \cite[Subsection 1.3]{AGR} for the study of
connected components of random graphs.

\textit{Proof of Proposition \ref{Skoro}}.
Set $\gamma_0=\tilde{\gamma}_0=0$, and, for $k\ge 1$:
\be
\tilde{\gamma}_k&=&\#\{\textrm{empty places in the set }
\{V,V+1,\ ...\ ,V+k-1\}\}\\
&=&\#\{j\ |\ V\leq j\leq V+k-1 \textrm{ and }H_j=0\}\\
&=&1+\#\{j\ |\ V+1\leq j\leq V+k-1 \textrm{ and }H_j=0\}\\
\gamma_k&=&\max_{1\leq i\leq k}\{C_{V-1}-
C_{V+i-1}\}.
\ee
Relation  (\ref{reccirc}) yields at once:
\be
H_{V+k}&=&H_{V}+Y_{V+1}
+Y_{V+2}+\ ...\ +Y_{V+k}-k+\tilde{\gamma}_k\\
&=&C_{V+k}-C_{V}+\tilde{\gamma}_k,
\ee
so the proof of Proposition \ref{Skoro} reduces to that  of $\tilde{\gamma}_k
=\gamma_k$.

We already have $\gamma_1=\tilde{\gamma}_1=1$.
For $k\ge 1$, note that either  $\gamma_{k+1}=\gamma_k+1$ or $\gamma_{k+1}=\gamma_k$.
First consider the case $\gamma_{k+1}=\gamma_k$: there exists $j$ such that
$0\le j\le k-1$ and $C_{V+j}\le C_{V+k}$.
This  can be rewritten:
\[Y_{V+j+1}+\ ...\ +Y_{V+k}\ge k-j,\]
meaning  that  more than $k-j-1$ cars want to park on only
$k-j$ places. Thus the last place, $V+k$, is necessarily occupied,
 i.e. $\tilde{\gamma}_{k+1}=\tilde{\gamma}_k$.

Assume now that $\gamma_{k+1}
=\gamma_k+1$: for any $j$ such that
$0\le j\le k-1$, we have $C_{V+j}\ge C_{V+k}+1$,
or equivalently:
\[Y_{V+j+1}+\ ...\ +Y_{V+k}\le k-j-1.\]
Using this inequality in the same way as we used relation (\ref{plein}) previously,
we conclude that if $\tilde{\gamma}_{k+1}
\not=\tilde{\gamma}_k+1$ or, equivalently, if $V+k$ is not empty, then
the sets of places  $\{V+j,\cdots,V+k\}$ have to  be full,
for any $j$ such that $0\le j\le k-1$, including thus $V$.
$\hspace{1cm}{\diamondsuit}$

We just proved that
\begin{pro}
Place $V+k$ is  empty if and only if
 $\gamma_{k+1}=\gamma_k+1$, or if and only if
$-C_j$ has a record at $j=V+k$.
\label{record}
\end{pro}
Sequence $\gamma_k$ will be easier to handle than
$\tilde{\gamma}_k$, when dealing with uniform convergence
 in the next subsection.

\subsection{Proof of Theorem \ref{monticule}}
\label{exc}

Recall that Donsker (1952),
 following an idea of Doob,
proved that:
\begin{theo}
\label{Donsk}
Let $b=(b(t))_{0\leq t\leq 1}$ be a Brownian bridge. We have:
\[\alpha_n\build{\tend}{}{weakly}b.\]
\end{theo}
We shall also need:
\begin{theo}(Vervaat, 1979 \cite{VER})
\label{waat}
Let $v$ be the almost surely unique point such that
$b(v)=\min_{0\leq t\leq 1}\ b(t)$.
Then $v$ is uniform and $e=(e(t))_{0\leq t\leq 1}$,
 defined by  $e(t)=b(\{v+t\})-b(v)$, is a
normalized Brownian excursion, independent of $v$.
\end{theo}

Owing to the Skorohod representation theorem
 \cite[II.86.1]{ROW}, we  assume the joint existence,
on some probabilistic triplet $(\Omega,A,P)$,
 of a sequence of copies of empirical processes, also denoted $\alpha_m$,
 and of a Brownian bridge $b$, such that,
for almost any $\omega \in \Omega$, $t\rightarrow\alpha_m(t,\omega)$
 converges uniformly on $[0,1]$ to $t\rightarrow b(t,\omega)$. We also assume,
in the definition of $h_{\lambda}$,
 that $e$ and $v$ are generated from $b$, using Vervaat's Theorem, so that
$h_{\lambda}=\Psi_{\lambda}b$.

The idea of the proof is to build a sequence of copies
of $h_m$ that converges almost surely uniformly to a copy of $h_{\lambda}$:
first  $\alpha_n$ defines sequences
\be
A^{m,\ell}_k&=&\sqrt n\hspace{0.1cm}\alpha_n(k/m),\\
C^{m,\ell}_k&=&\sqrt n\hspace{0.1cm}\alpha_n(k/m)
-\ell\hspace{0.1cm}\frac{k}{m}.
\ee
Then, from $C^{m,\ell}_k$, we can  define, through Proposition \ref{Skoro},
$H^{m,\ell}_k$ that is distributed as
$H_k$, though no underlying parking scheme has been defined.
Actually we can also define
\[\left((A^{m,\ell}_k,C^{m,\ell}_k,H^{m,\ell}_k,\gamma^{m,\ell}_k,
Y^{m,\ell}_k)_{-\infty\le k\le +\infty},V(m,\ell)\right),\]
with the same distribution as
$\left((A_k,C_k,H_k,\gamma_k,Y_k)_{-\infty\le k\le +\infty},V\right)$
 in the previous Subsection.
Thus $\tilde h_m$ defined by
\[\tilde h_m(t)=\frac{1}{\sqrt m}H^{m,\ell}_{\lfloor mt\rfloor}\]
is distributed as $h_m$, and we shall drop the tilda in
 what follows.  We also set:
\begin{eqnarray}
\nonumber z_m(t)&=&\frac{1}{\sqrt m}\left(C^{m,\ell}_{V(m,\ell)+\lfloor mt\rfloor}-
C^{m,\ell}_{V(m,\ell)}\right).
\end{eqnarray}

We have
\begin{lem}
\label{fourretout}
If
$\lim_{m} \frac{\ell}{\sqrt m} = \lambda,$
then for almost any $\omega$,
\begin{eqnarray}
\label{stepconstant}
\alpha_n(\lfloor mt\rfloor/m)\hspace{0,3cm}
&\build{\tend}{}{uniformly}&
\hspace{0,3cm} b(t);\\
\label{min}
\lim_m\frac{V(m,\ell)}{m}&=&v;\\
z_m(t)\hspace{0,3cm}&\build{\tend}{}{uniformly}&
\hspace{0,3cm}z(t)=e(t)-\lambda t;
\label{excur}\\
h_m(\{t+(V(m,\ell)/m)\})\hspace{0,3cm}
&\build{\tend}{}{uniformly}&
\hspace{0,3cm}\Psi_{\lambda}e(t).
\label{visites}
\end{eqnarray}
\end{lem}

Theorem \ref{monticule} is a
reformulation of
(\ref{visites}), since we have:
\[\Psi_{\lambda}e(\{t-(V(m,\ell)/m)\})\hspace{0,3cm}
\build{\tend}{}{uniformly}
\hspace{0,3cm}\Psi_{\lambda}e(\{t-v\}) = h_{\lambda}(t).\]

\textit{Proof of  (\ref{stepconstant})}. Set
 \[M_m=\max_{0<k\leq m}Y^{m,\ell}_k.\]
We have:
\be
|\alpha_n(\lfloor mt\rfloor/m)-\alpha_n(t)|&\leq&
\frac{\sqrt n}{m}+\frac{M_m}{\sqrt n}\\
&\leq&
\frac{1+M_m}{\sqrt n},
\ee
and, as $Y^{m,\ell}_k$ follows the  binomial distribution
 with parameters $(n,1/m)$,
\begin{eqnarray}
\nonumber
\Pr(M_m\geq C\log m)&\leq&m\Pr(Y^{m,\ell}_1\geq C\log m)\\
\nonumber&\leq&mE[\exp(KY^{m,\ell}_1)]\exp(-KC\log m)\\
\label{variation}&\leq&Am^{1-KC}.
\end{eqnarray}
Thus Borel-Cantelli Lemma entails
that, for a suitable $C$, with probability $1$
the supremum norm of
$\alpha_n(\lfloor mt\rfloor/m)-\alpha_n(t)$
vanishes as quickly as $\frac{C\log m}{\sqrt n}.
\hspace{1cm}\diamondsuit$

\textit{Proof of  (\ref{min})}. For this proof
 and the next one,
we consider an $\omega$ such that simultaneously
$\alpha_n(t,\omega)$ and $\alpha_n(\lfloor mt\rfloor/m,\omega)$
converges uniformly (for $t\in [0,1]$) to $b(t,\omega)$, and
such that $t\rightarrow b(t,\omega)$ reaches its minimum
 only once (we know that the set of such $\omega$'s has
 measure $1$). We set:
\be
\varepsilon_{m,1}&=&\sup_{0\leq t\leq 1}
|\alpha_n(\lfloor mt\rfloor/m)-\alpha_n(t)|,\\
\varepsilon_{m,2}&=&\sup_{0\leq t\leq 1}|b(t)-\alpha_n(t)|,\\
\varepsilon_{m,3}&=&\sup_{0\leq t\leq 1}\left\vert b(v+t)-
b\left(\frac{V(m,\ell)}{m}+t\right)\right\vert.
\ee

From the continuity property of $b$,
 the first minimum, $V(m,\ell)/m$, of $\alpha_n(\lfloor mt\rfloor/m)$
 converges to the only minimum of $b$
 (i.e. $v$): we have
\be
b(V(m,\ell)/m)&\leq&\alpha_n(V(m,\ell)/m)+\varepsilon_{m,2}\\
&\leq&\alpha_n(\lfloor mv\rfloor/m)+\varepsilon_{m,2}\\
&\leq&b(v)+\varepsilon_{m,1}+2\varepsilon_{m,2}.
\ee
Now the minimum of $b(t)$ over the set $[v-\varepsilon,
v+\varepsilon]^c\cap[0,1]$ is $b(v)+\eta$ for some positive $\eta$,
and thus, if $\varepsilon_{m,1}+2\varepsilon_{m,2}<\eta$, then necessarily
$|v-V(m,\ell)/m|<\varepsilon$
.$\hspace{1cm}\diamondsuit$

\textit{Proof of (\ref{excur})}. One checks easily that:
\[
|z_m(t)-z(t)|\leq 2(\varepsilon_{m,1}+\varepsilon_{m,2}+\varepsilon_{m,3})+\frac{\lambda }{m}
+\left\vert\lambda-\frac{\ell}{\sqrt n}\right\vert.\hspace{1cm}\diamondsuit\]

\textit{Proof of (\ref{visites})}.
 According to  Proposition \ref{Skoro},
we have:
\be
h_m\left(\{t+(V(m,\ell)/m)\}\right)&=&z_m(t)+
\max_{0\leq s\leq t}\left(-z_m(s)\right)+\varepsilon_{m,4}(t),
\ee
where
\be
\varepsilon_{m,4}(t)&=&\max_{0\leq s\leq t}\left(
\frac{Y^{m,\ell}_{V(m,\ell)+\lfloor ms\rfloor}}{\sqrt n}
-z_m(s)\right)-\max_{0\leq s\leq t}(-z_m(s)),\\
|\varepsilon_{m,4}(t)|&\leq&\frac{M_m}{\sqrt n}.
\ee
Thus  (\ref{visites})
 follows from the uniform convergence
of $z_m$ to $z$.$\hspace{1cm}\diamondsuit$

\subsection{An extension of the invariance principle}
\label{ald97}

This section is the last step of the proof of Theorem \ref{sorted}.
 The widths of excursions of $h_m(t)$
above zero are the sizes of the blocks of cars of
 the corresponding parking scheme, normalized by $m$.
Unfortunately, uniform convergence of $h_m$ to $h$
 does not entails  convergence of sizes of excursions.
However the excursions of $z_m$ above its current minimum are
 exactly the excursions of $h_m$ above $0$, up to the random
shift $V(m,\ell)/m$, and, according to \cite[Section 2.3]{AGR},
the uniform convergence of $z_m$ to $z$ entails convergence of
 sizes of excursions of $z_m$ above its current minimum
to  sizes of excursions of $z$ above its current minimum,
 provided that $z$ does never reach
its current minimum \textit{two times}.
It is known that this last condition holds true
for almost each sample path $z$,
so that we have almost sure convergence
of sizes of excursions of $z_m$, or equivalently of sizes of blocks.
Similarly, excursions of $z$ above its current minimum are
also excursions of $\Psi_{\lambda}e$
 above $0$,  yielding Theorem \ref{sorted}.

Let us give some details and notations.
We shall apply to $z_m$ and $z$
the following weakened form of
\cite[Lemma 7, p. 824]{AGR}:

\begin{lem} Suppose $\zeta$:
$[0,+\infty[\tend R$ is continuous.
 Let $E$ be the set of nonempty
 intervals $I=(l,r)$ such that:
\[\zeta(r)=\zeta(l)= \min_{s\leq l}\zeta(s),
\hspace{2cm} \zeta(s)>\zeta(l)
\textrm{ for }l<s<r.\]
Suppose that, for intervals
 $I_1$, $I_2\in E$ with
$l_1<l_2$ we have \[\zeta(l_1)>\zeta(l_2).\]
Suppose also that the complement of
$\cup_{I\in E}(l,r)$ has
 Lebesgue measure $0$. Let
$\Theta = \{(l,r-l):\ (l,r)\in E\}$.
 Now let $\zeta_m\tend \zeta$
 uniformly on $[0,1]$.
 Suppose $(t_{m,i}$, $i\geq 1)$
satisfy the following:
\begin{itemize}
\item[($i$)] $0 = t_{m,1}<t_{m,2}
< ... <t_{m,k+1}=1;$
\item[($ii$)] $\zeta_m(t_{m,i})=\min_{u\leq t_{m,i}} \zeta_m(u);$
\item[($iii$)] $\lim_m\max_i(\zeta_m(t_{m,i})-\zeta_m(t_{m,i+1}))=0$
\end{itemize}
Write $\Theta^{(m)}=\{(t_{m,i},t_{m,i+1}-t_{m,i});
\ 1\leq i\leq k\}$. Then $\Theta^{(m)}\tend \Theta$
for the vague topology of measures on
$[0,1]\times (0,1]$.
\label{lem7}
\end{lem}

Set of points, such as $\Theta^{(m)}$ or $\Theta$, can also be seen as
 \textit{point processes} (i.e. measures that are  infinite sums of Dirac masses): we
 identify the set $A$ and the measure
\[\sum_{x\in A}\delta_x.\]
For point processes on
$[0,1]\times (0,1]$, the following criterium of convergence holds:
\begin{pro}
$\Theta^{(m)}\tend \Theta$
 for the vague topology if and only if,
 for any $y>0$ such that $\Theta([0,1]\times\{y\})=0$,
\begin{itemize}
\item[($i$)]
 for $m$ large enough, $\Theta^{(m)}([0,1]\times [y,1])
=\Theta([0,1]\times [y,1])$;
\item[($ii$)] for any $x\in[0,1]\times [y,1]$
such that $\Theta(\{x\})>0$ there
 is a sequence of points $x_m$,
$\Theta^{(m)}(x_m)>0$, such that  $x_m\tend x$.
\end{itemize}
\label{vague}
\end{pro}

As an easy consequence, partly due to
 the fact that second components add up to $1$:
\begin{cor}
If $\Theta^{(m)}\tend \Theta$
 for the vague topology, then the sequence of
second components of points of $\Theta^{(m)}$,
sorted in decreasing order,
converge componentwise and in
$\ell_1$ to the corresponding sequence for $\Theta$.
\label{ranked}
\end{cor}

One can find the proofs of Lemmata and Propositions of this
subsection, and also of the stochastic calculus
points in the next proof, in \cite[pp. 30-34]{CL}.

Let us choose $(\zeta_m,\zeta)=(z_m,z)$,
defined at Subsection \ref{exc}.
Let the $t_{m,i}$'s of Lemma \ref{lem7}
 be the successive positive  records of $-z_m(t)$
so, due to Lemma \ref{record}, the
$V(m,\ell)+mt_{m,i}$'s are the  $\ell$ empty places
of the corresponding parking scheme,
 counted starting at $V(m,\ell)$.
The sequence of
second components of $\Theta^{(m)}$ (resp. of $\Theta$)
is nothing else but $\frac{1}{m}B^{m,\ell}$
(resp. $B(\lambda)$).
Thus Theorem  \ref{sorted} follows from
Lemma \ref{lem7} and Corollary \ref{ranked},
applied to $\zeta(t)=z(t)$,
$\zeta_m(t)=z_m(t)$.
Let us check the assumptions
of Lemma \ref{lem7}. First, not depending on $i$,
$$z_m(t_{m,i})-z_m(t_{m,i+1})=1/\sqrt m,$$
giving  assumption $(iii)$.
The standard Brownian motion satisfies the assumption
 "almost surely, $\zeta(l_1)<\zeta(l_2)$ for any $l_1<l_2$", and,
due the Cameron-Martin-Girsanov formula,  this extends to
solutions of  stochastic differential equations with smooth coefficients,
including $z$ (cf.  \cite[Chp. XI, Ex. 3.11]{RY}).
Setting $O=\cup_{I\in E}(l,r)$,
the Lebesgue measure of $O^c$ is $0$
for similar reasons
(see \cite[pp. 33-34]{CL} for details).
$\hspace{1cm}\diamondsuit$

\section{Extension to finite-dimensional distributions}
\label{uniform}

This Section is devoted to the proof of Theorem \ref{lambeaux}.
 Up to now, with the exception of Subsection \ref{coal}, we only considered
the parking process frozen at a given time $n$, that is, just after the arrival of car $c_n$.
Theorem \ref{lambeaux} is a result about the dependence between parking schemes,
at successive times $n_1<n_2<\cdots<n_k$. Thus we shall
need a two-parameters (time and place)  analog of Theorem \ref{monticule}.
For each $m,\ell$, let $h_{m,\ell}(t)$ be  the profile
of  the parking scheme
 of the $m-\ell$ first cars on the
$m$ places. Similarly,  let  $z_{m,\ell}(t)$
be the analog of $z_m$ defined at Section \ref{profil}.
Finally, for $\lambda\sqrt m\le m$, set
\[\psi_m(\lambda,t)=h_{m,\lceil \lambda\sqrt m\rceil}(t),\]
else let $\psi_m(\lambda,t)=0$.
The  dependence between the $m$ successive parking schemes, after the $m$ successive
arrivals on $m$ places is captured by the two-parameters process
\[\psi_m=\left(\psi_m(\lambda,t)\right)_{0\le\lambda,\,0\le t\le 1}.\]
\begin{figure}[ht]
\vspace{-0.5cm}
\hspace{-0.5cm}
\includegraphics[width=14cm]{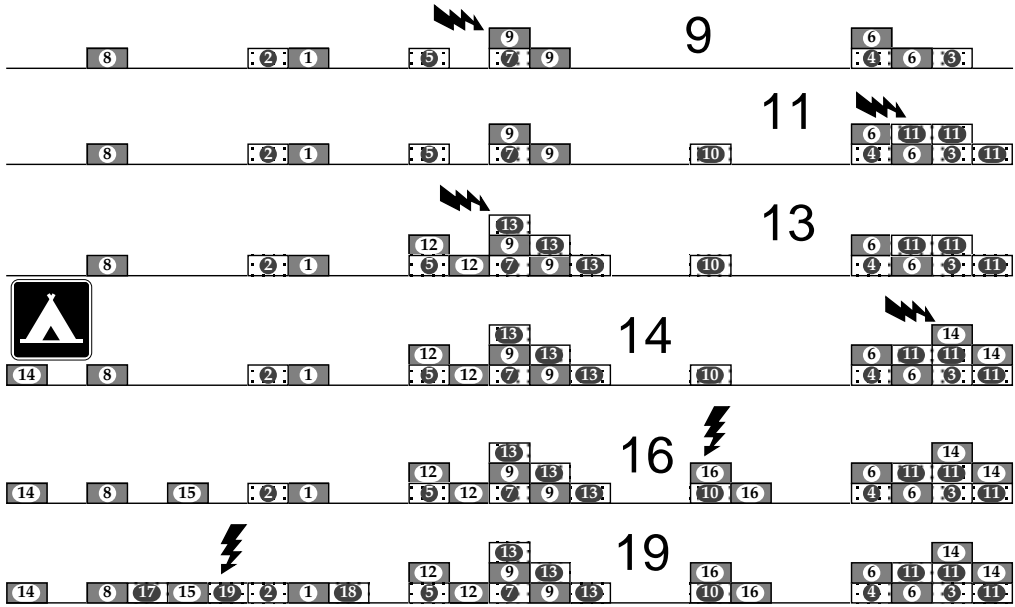}
\caption{Parking schemes for  $m=25$ places
 and $n=1, \cdots,19$ cars.}
\end{figure}
Note that the time parameter, $\lambda$, decreases as time goes by and cars
arrive, while $t$ is the location parameter: $\sqrt m\ \psi_m\left(\frac{\ell}{\sqrt m},
\frac{k}{m}\right)$ is the number of cars that tried to park, successfully or not,
 on place $k$, among the $m-\ell$ cars already arrived. We have:
\begin{theo}
\label{2dim}
There exists, on some probability space $\Omega$, a uniform random variable
 $v$, and copies of  $\psi_m$ and of the normalized Brownian excursion $e$, such that,
for $D_{\Lambda}=[0,\Lambda]\times[0,1]$,
\[\Pr\left(\forall \Lambda,\hspace{0.2cm}\psi_m(\lambda,t)\hspace{0,3cm}
\build{\tend}{on\textrm{ }D_{\Lambda}}{uniformly}
\hspace{0,3cm}h_{\lambda}(t)\right)=1.\]
\end{theo}
Set
\[Z_m(\lambda,t)=z_{m,\lceil \lambda\sqrt m\rceil}(t).\]
We shall actually prove that
\begin{equation}
\label{Z}
\Pr\left(\forall \Lambda,\hspace{0.2cm}Z_m(\lambda,t)\hspace{0,3cm}
\build{\tend}{on\textrm{ }D_{\Lambda}}{uniformly}
\hspace{0,3cm}e(t)-\lambda t\right)=1.
\end{equation}
Theorem \ref{2dim} will follow, as well as a description
of the asymptotic  evolution, as cars arrive, of the whole sequence of
sizes \textit{and positions} of blocks.

We need more notations
to give a precise statement. Let $\Theta^{(m)}(\lambda)$
 denote the point process
corresponding to the choice
 $\zeta_m=Z_m(\lambda,\cdot)$
in Lemma \ref{lem7}:
the first components of  points of $\Theta^{(m)}(\lambda)$
are the positions, \textit{relative to} $V(m,\ell)$ and normalized by $m$,
of the $ \lceil \lambda\sqrt m\rceil$ empty places after the
 $ \lfloor m-\lambda\sqrt m\rfloor^{th}$ arrival; the second components
are the lengths, normalized by $m$,
 of blocks starting at these empty places (the length including also
the initial empty place). We allow empty blocks, that is, empty places
followed by another empty place: the corresponding length is $1/m$.
Similarly, let $\Theta(\lambda)$ denote the point process
corresponding to the choice
 $\zeta(t)=e(t)-\lambda t$
in Lemma \ref{lem7}: the first component of an element
 of $\Theta(\lambda)$ is the starting point of an excursion of $\Psi_
{\lambda}e$,  the second component of this element being
 the width of the same excursion. We have
\begin{theo}
\label{totale}
The finite-dimensional distributions of $\Theta^{(m)}(\lambda)$
converges weakly to the finite-dimensional distributions of $\Theta(\lambda)$.
\end{theo}

This result is weaker than the weak convergence of  $\Theta^{(m)}$
to $\Theta$, that is, it does not insure the weak convergence of any continuous
functional of $\Theta^{(m)}$ to the same functional applied to $\Theta$,
but it insures that if $\Theta^{(m)}$ has a weak limit, this limit can only be
 $\Theta$.

\textit{Proof of Theorem \ref{2dim}.}
As in Subsection \ref{exc}, we start, on some space $\Omega$,
with a sequence $\alpha_m$ of empirical processes that
converges almost surely uniformly to a Brownian  bridge $b(t)=e(\{t-v\})-
e(-v)$. For the proof of  Theorem \ref{monticule}, there was no need
to  build a  random parking scheme corresponding to $\alpha_m$ - but, maybe,
 for the mental picture. This task cannot be avoided now, as we need the
 \textit{chronology} to deduce $h_{m,\ell}$, $z_{m,\ell}$,
$\psi_m$, $Z_m$ and $\Theta_m$ from $\alpha_m$.

There is however a slight difficulty: $\alpha_m$ provides
the total number $Y^{m,\ell}_k$ of cars whose first
try was on place $k$, but it does not provide the chronology.
Let us collect some basic facts concerning empirical processes:
 $\alpha_m$ has $m$ positive jumps with height $\frac{1}{\sqrt m}$,
 at places that we call $\left(J_k^{(m)}\right)_{1\leq k\leq m}$.
Between the jumps $\alpha_m$ has the  slope - negative  -
$-\sqrt m$. The random vector $J^{(m)}=\left(J_k^{(m)}\right)_{1\leq k\leq m}$
is uniformly distributed on the simplex
$\{0<x_1<x_2< ... <x_m<1\}$. Any random permutation $\sigma_m$
of  $J^{(m)}$'s components,  with $\sigma_m$ and $J^{(m)}$ independent, yields
a sequence
\[\left(U_k^{(m)}\right)_{1\leq k\leq m}
=\left(J_{\sigma_m(k)}^{(m)}\right)_{1\leq k\leq m}\]
of independent uniform
random variables on $[0,1]$, whose empirical process is $\alpha_n$.

Thus we  can
recover the chronology with the help of $\sigma_m$,
assuming that car $c_k$
tries to park first on place $\left\lceil nU_k^{(m)}\right\rceil$.
Let us define $\alpha_{m,\ell}(t)$ (resp. $\tilde\alpha_{m,\ell}(t)$)
as the empirical processes
for the samples $(U^{(m)}_i)_{1\le i\le m-\ell}$
(resp. $(U^{(m)}_i)_{m-\ell+1\le i\le m}$).
Both are samples of
independent and uniform random variables. Now $\alpha_{m,\ell}(t)$ allows to define
the profiles $h_{m,\ell}(t)$ of the successive parking schemes, and to define 
$z_{m,\ell}(t)$, $\psi_m$  and $Z_m$ as well, following the same lines as
in Subsection \ref{exc}. 

We shall see now that  any choice of the
sequence $(\sigma_m)_{m\ge 1}$ of uniform random permutations
 insures the convergence of $\Theta_m$ to $\Theta$, provided that  $\alpha_m$
and $\sigma_m$ are independent for each $m$. 
We give at the end of the proof a  construction of $\sigma_m$ that 
 will be useful in Section \ref{fractal}. 
We have:
\[\sqrt{m}\hspace{0,2cm}\alpha_m(t)=\sqrt{m-\ell}
\hspace{0,2cm}\alpha_{m,\ell}(t)+\sqrt{\ell}
\hspace{0,2cm}\tilde\alpha_{m,\ell}(t),
\]
with the consequence that:
\[|\alpha_m(t)-\alpha_{m,\ell}(t)|\le
\left|-1+\sqrt{1-\frac{\ell}{m}}\right|
|\alpha_{m,\ell}(t)|+\sqrt{\frac{\ell}{m}}
\hspace{0,1cm}|\tilde\alpha_{m,\ell}(t)|.
\]
According to the DKW inequality \cite{MASS},
not depending on $(m,\ell)$,
\[\Pr(\sup_t|\tilde\alpha_{m,\ell}(t)|\geq x)\leq
2 \exp(-2x^2),\]
thus,
for suitable $K_1$ and $K_2$, and for $\varepsilon >0$,
\[\Pr\left(\sup_{0\leq \ell\leq \Lambda\sqrt m}
\sup_t|\alpha_m(t)-\alpha_{m,\ell}(t)|
\ge m^{-1/4+\varepsilon}\right)\le
K_1\sqrt m\ e^{-K_2m^{2\varepsilon}}.\]
Thus, using Borel-Cantelli lemma, we obtain
\begin{equation}
\label{unifff}
\Pr\left(\sup_{0\leq \ell\leq \Lambda\sqrt m}
\sup_t|\alpha_m(t)-\alpha_{m,\ell}(t)|
=O(m^{-1/4+\varepsilon})\right)=1.
\end{equation}

Owing to (\ref{unifff}), a simple glance at
the proof of  (\ref{min}) show that
the convergence
 of  $V(m,\ell)/m$ to $v$,  for $0\leq \ell\leq \Lambda\sqrt m$,
  is uniform, almost surely. Slightly changing the
 definitions of $\varepsilon_{m,1}$, $\varepsilon_{m,2}$, $\varepsilon_{m,3}$ of Subsection
\ref{exc}  and defining also $\varepsilon_{m,5}$, $\varepsilon_{m,6}$, $v_m(\lambda)$
 as follows:
\be
\varepsilon_{m,1}&=&\sup_{0\leq t\leq 1}
|\alpha_m(\lfloor mt\rfloor/m)-\alpha_m(t)|,\\
\varepsilon_{m,2}&=&\sup_{0\leq t\leq 1}
|b(t)-\alpha_{m}(t)|,\\
\varepsilon_{m,3}&=&\sup_{0\leq \ell\leq
\Lambda\sqrt m}\ \ \sup_{0\leq t\leq 1}|b(t+v)-b(t+V(m,\ell)/m)|,\\
\varepsilon_{m,5}&=&\sup_{0\leq \ell\leq
\Lambda\sqrt m}\  \ \sup_{0\leq t\leq 1}
|\alpha_{m}(t)-\alpha_{m,\ell}(t)|,\\
\varepsilon_{m,6}&=&\sup_{0\leq \lambda\leq\Lambda}\ \ \sup_{0\leq t\leq 1}
\left|\frac{\lceil \lambda\sqrt m\rceil}
{\sqrt {m-\lceil \lambda\sqrt m\rceil}}
\frac{\lfloor mt\rfloor}{m}-\lambda t\right|\\
v_m(\lambda)&=&\frac{V(m,\lceil \lambda\sqrt m\rceil)}{n},
\ee
we have, for $i\in\{1,2,3,5,6\}$,
\[\lim_m\ \varepsilon_{m,i}=0,\]
from Subsection
\ref{exc} for $i=1,2$ and from (\ref{unifff}) and uniform continuity of $b$
for $i=3$. Furthermore, we have
\be
|Z_{m}(\lambda,t)-e(t)+\lambda t|&\le&\varepsilon_{m,1}+\varepsilon_{m,2}+
\varepsilon_{m,5}+\varepsilon_{m,6}+\left|b\left(t+v_m(\lambda)\right)
-b(t+v)\right|\\ &\le&\varepsilon_{m,1}+\varepsilon_{m,2}+\varepsilon_{m,3}+
\varepsilon_{m,5}+\varepsilon_{m,6}.
\ee

Finally let us give  a  construction of $\sigma_m$ that 
 will prove useful in Section \ref{fractal}. We  can enlarge 
 the probability space $\Omega$,
provided by the Skorohod representation Theorem, to
$\Omega\times[0,1]^{\{1,2,3,\dots\}}$, obtaining a  sequence of independent
random variables $(u_k)_{k\geq 1}$, uniform  on $[0,1]$ and independent
 of the sequence $(\alpha_m)_{m\geq 1}$, and we let
\[\sigma_m(k)=\#\left\{i\ \vert\  1\le i\le m\textrm{ and }u_i\le u_k\right\}.\]
Note that with this construction of  the sequence $U^{(n)}$,  $U^{(n)}$ cannot
 be obtained by erasing the last term of $U^{(n+1)}$,
as usual. 
If it was the case,
$\alpha_n$ would not converge uniformly to $b$, due to
 Finkelstein's law of the iterated logarithm \cite[Theorem 5.1.2]{CSO}.
Incidentally, Finkelstein's law suggests that $\alpha_{m,\ell}(t)$  converges
uniformly to $b$ only if we choose $\ell=o(m)$.$\hspace{1cm}\diamondsuit$

\textit{Proof of Theorem \ref{totale}.}
We have seen that, almost surely, for a given $\lambda$,
the assumptions of Lemma \ref{lem7} hold true for $\zeta(t)=
 e(t)-\lambda t$, so, still almost surely, they hold true jointly for
 $0\le\lambda_1\le\lambda_2\le\cdots\le\lambda_k$, yielding that
\[\Pr\left(\lim_m\left(\Theta^{(m)}(\lambda_i)\right)_{1\le i\le k}
\hspace{0,3cm}=\hspace{0,3cm}\left(\Theta(\lambda_i)\right)_
{1\le i\le k}\right)=1.\]

\textit{Proof of Theorem \ref{lambeaux}}.
\[\lim_m\left(\Theta^{(m)}(\lambda_i)\right)_
{1\le i\le k}\hspace{0,3cm}
=\hspace{0,3cm}\left(\Theta(\lambda_i)\right)_
{1\le i\le k}\]
entails that
\[\lim_m\left(B^{(m)}(\lambda_i)\right)_
{1\le i\le k}\hspace{0,3cm}
=\hspace{0,3cm}\left(B(\lambda_i)\right)_
{1\le i\le k}.\]

\section{Distribution of components of $B(\lambda)$ and $R(\lambda)$}
\label{Markov}
The proofs of Theorems \ref{Pavlov}  and  \ref{size-biased}, that we give in this
 Section, are more  of a combinatorial nature.

\subsection{Proof of Theorem \ref{Pavlov}}
\label{prPavlov}

This proof reduces to explain a  one-to-one correspondence between
\textit{confined} parking schemes  with $n$ cars and $\ell$ empty places
and Pavlov's forests with  $\ell$ rooted trees and $n$ non-root vertices,
correspondence in which the sizes of trees and the sizes
  of blocks are in correspondence too.
Then  Theorem  \ref{Pavlov}
is just a restatement of \cite[Theorem 4]{PAV}.

\begin{figure}[ht]
\vspace{-2.2cm}
\hspace{-1cm}
\includegraphics[width=14cm]{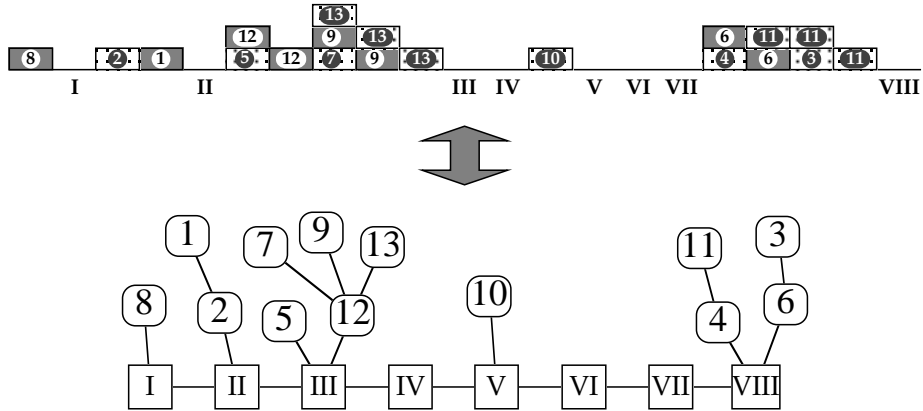}
\vspace{0,4cm}
\caption{Correspondence parking $\leftrightarrow$
 Pavlov's forests,
empty places $\leftrightarrow$ roots.}
\label{forest}
\end{figure}
In Pavlov forests, roots and  non-roots are  labeled separately, for instance
the roots  (resp. non-roots) are labeled $\{r_1,r_2,\ \dots\ ,r_{\ell}\}$ (resp.
$\{v_1,v_2,\ \dots\ ,v_n\}$). The label of  the root is also the label of the corresponding tree.
Let us define the Pavlov forest $T$ corresponding  to a given confined parking
scheme $\Pi$: the non-roots of the first tree of $T$ (that is, of the tree  rooted at $r_1$)
are the cars parked before the first empty place of $\Pi$, and the non-root vertices of the
$k^{th}$  tree are the $s_k$ cars parked between  the $k-1^{th}$  and the $k^{th}$
empty places. The way these $s_k$ cars are parked  can be described by a
confined parking scheme of $s_k$ cars on $s_k+1$ places: we define the  $k^{th}$  tree
of $T$ through one among the many one-to-one
correspondences  between rooted labeled trees with
$m$ nodes, and confined parking schemes of $m-1$ cars on $m$
 places  \cite{CM,FOA,FRA,SHU}.

The following one-to-one correspondence will be specially useful at Section \ref{concl},
to explain the relation between  parking and the standard additive coalescent. Consider
 a random labeled tree $t$ with $k$  vertices  $v_1,\ \dots\ ,v_k$ and let $\Pi$ denote the
 corresponding  confined parking scheme for $k-1$ cars
\begin{figure}[ht]
\begin{center}
\includegraphics[width=7.5cm]{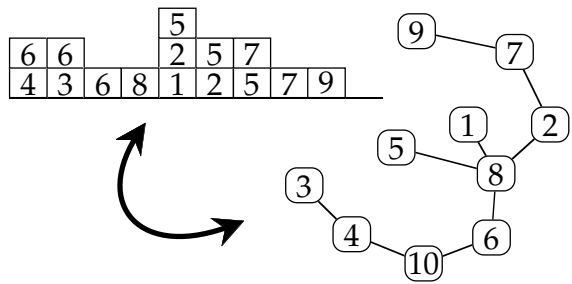}
\end{center}
\caption{Correspondence parking $\leftrightarrow$
 labeled tree}
\label{bijection-arbres}
\end{figure}
 on $k$ places. The description of $\Pi$ uses a variant of the breadth first search of $t_k$:
by convention $v_k$ is the root of $t$  ; at step $1$, $v_k$'s sons  are stored
in a queue, \textit{the smallest labels at the head of the queue}. Then at each step the vertex
 at the head of the queue is removed from the queue, while its sons are added to the queue,
and \textit{the queue is reordered} (the smallest labels at the head) to be ready for the
 following step. The corresponding parking scheme  $\Pi$ is defined by specifying that the
 first try of car $c_j$ is on place $i$ if and only if the first appearance of $v_j$ in the queue
 is at step $i$. In this correspondence, one checks easily that car $c_j$ finally parks at place $i$ if and only
 if $v_j$ is at the head of the queue at step $i$,  and also that the successive lengths of the  queue
 just give $\Pi$'s profile.

\subsection{Proof of Theorem \ref{size-biased}}
\label{prsize-biased}
We first provide a useful identity
leading to the proof of Theorem \ref{size-biased}. Set
\[d_i=k_1+k_2+ ... +k_i.\]
We have
\begin{pro}
\label{multiple}
\[\Pr(R_j^{m,\ell}=k_j,\hspace{0,2cm}
1\leq j\leq i)=\prod_{j=0}^{i-1}\varphi(m-d_j-j,n-d_j,k_{j+1}).\]
\end{pro}
\pr The choice of the elements in each of the blocks can be done
in
\[\prod_{j=1}^{i}C^{k_j-1}_{n-d_{j-1}-1}\]
ways, and they can be arranged inside each of these blocks in
\[\prod_{j=1}^{i}(k_j+1)^{k_j-1}\]
ways.

It will be  convenient to argue in terms of confined
 parking schemes, since rotations do not change the sizes of blocks.
The total number of confined parking schemes
is $m^{n-1}\ell$.
We obtain a confined parking scheme with
blocks' sizes $k_1$, $k_2$, etc ... ,
for the $i$ first blocks, respectively,
by inserting these $i$ blocks successively, with an empty
place attached to the right of them, insertion taking place at the
front of the confined parking scheme $\Pi$ for  the remaining cars, or
just after one of the empty
places of the confined parking scheme for the remaining cars.
There are $(m-d_i-i)^{n-d_i-1}(m-n-i)$ choices for $\Pi$,
$m-n-i+1$ possible insertions for the first block, $m-n-i+2$
possible insertions
for the second block, and so on ...
Finally, the probability
$p(k)$ on the left hand of Proposition \ref{multiple}
is given by
\[p(k)=\frac{(m-d_i-i)^{n-d_i-1}
(m-n-i)}{m^{n-1}(m-n)}
\prod_{j=1}^{i}(k_j+1)^{k_j-1}(m-n+j-i)
C^{k_j-1}_{n-d_{j-1}-1}.\]
It is not hard to check
 that this last expression
 is the same as the right
 hand of Proposition \ref{multiple}.
$\hspace{1cm}\diamondsuit$

\textit{Proof of Theorem \ref{size-biased}.} Set $S_0=s_0=0$ and
\be
S_j&=&R_1(\lambda)+R_2(\lambda)+ ... +R_{j}(\lambda),\\
s_i&=& x_1+x_2+ ... +x_i.
\ee
According to Theorem \ref{mixed} (that
will be proved independently at Subsection \ref{p1.5}), the
law of $(R_1(\lambda),R_2(\lambda), ... ,R_i(\lambda))$ has the following
 alternative characterization: for any $k$, conditionally, given $(R_1(\lambda),
R_2(\lambda), ... ,R_k(\lambda))$,   $R_{k+1}(\lambda)$ is distributed as
\[\left(1-S_k\right)\frac{N_{k+1}^2}{\frac{\lambda^2}{1-S_k}+N_{k+1}^2},\]
 in which  $N_{k+1}$ is standard Gaussian and independent of $S_k$. In other terms,
$R_{k+1}(\lambda)$  has the following conditional density:
\[\frac{1}{1-S_{k}}f\left(\frac{\lambda}{\sqrt{1-S_{k}}},
\frac{x_{k}}{1-S_{k}}\right),\]
and $R_{1}(\lambda)$  has the unconditional density $f(\lambda,x)$.
On the other hand, using the same line of proof as in Theorem \ref{undemi},
the approximations of Lemma \ref{local} for $\varphi(m,n,k)$ and
 Proposition \ref{multiple} lead,  for $\left(R_1^{(m)}(\lambda),R_2^{(m)}(\lambda),
 ... ,R_i^{(m)}(\lambda)\right)$, to  the following limit  density:
\[\prod_{j=0}^{i-1}\frac{1}{1-s_j}f\left(\frac{\lambda}{\sqrt{1-s_j}},
\frac{x_{j+1}}{1-s_j}\right).\hspace{1cm}\diamondsuit\]

\section{Sampling excursions of $\Psi_{\lambda}e$}
\label{fractal}

In this Section, we give the proofs of  Theorems \ref{fidis2} and \ref{decomp}.
They make essential use, to build the parking schemes,
of the  random permutation of  jumps of $\alpha_m$
defined at Section \ref{uniform}.

\subsection{Proof of Theorem \ref{fidis2}. }
\label{prfidis2}
We build a probability space
where almost sure convergence of $R^{(m)}_1(\lambda)$ to $R_1(\lambda)$ holds for
 any $\lambda$. As a consequence, for any $k$ and any $(\lambda_1,\lambda_2,
\ \dots\ ,\lambda_k)$,
\[\left(R^{(m)}_1(\lambda_1),R^{(m)}_1(\lambda_2),\ \dots\ ,R^{(m)}_1(\lambda_k)\right)
\build{\tend}{}{a.s.}{}\left(R_1(\lambda_1),R_1(\lambda_2),
\ \dots\ ,R_1(\lambda_k)\right),\]
entailing the result.

As in Section \ref{uniform},  we  enlarge  the probability space to
$\Omega\times[0,1]^{\{1,2,3,\dots\}}$, obtaining a  sequence of
i.i.d. random variables $(u_k)_{k\geq 1}$, uniform  on $[0,1]$ and independent
 of  $b$ and $(\alpha_m)_{m\geq 1}$. We let
\be
\sigma_m(k)&=&\#\left\{i\ \vert\  1\le i\le m\textrm{ and }u_i\le u_k\right\},\\
U_k^{(m)}&=&J_{\sigma_m(k)}^{(m)},\\
\pi_m&=&
\left\{U_1^{(m)}-v_m(\lambda)\right\},
\ee
and $\rho_1=\{u_1-v\}$, so that $\rho_1$ is uniform and independent of $e$.
We still assume that  first try of car $c_k$ is  place $\left\lceil mU_k^{(m)}
\right\rceil$. Thus, counted from $V\left(m,\left\lceil\lambda\sqrt m
\right\rceil\right)$, car $c_1$ parks at place $\left\lceil m\pi_m
\right\rceil$. Borel-Cantelli Lemma yields that
\begin{lem}
\label{mam dosz}
Almost surely,
\[\lim_m\ \left(U_1^{(m)},\pi_m\right)=(u_1,\rho_1).\]
\end{lem}

Let $g_m(\lambda)$ (resp. $d_m(\lambda)$) be the last zero
 of $t\rightarrow \tau_{v_m(\lambda)}\psi_m(\lambda,t)$ on the left of $\pi_m$
(resp. the first zero on the right). That is, $mg_m(\lambda)-1$
(resp. $md_m(\lambda)$)  is the empty place  at the  beginning (resp. at the end)
of the block containing car $c_1$, \textit{counted from} $V(m,\lceil\lambda\sqrt m
\rceil)$. Thus
\[d_m(\lambda)-g_m(\lambda)=R^{(m)}_1(\lambda).\]

Almost surely, due to Lemmata  \ref{fourretout} and  \ref{mam dosz},
\be
\lim_m\  Z_m(\lambda,\pi_m)&=&z(\rho_1).
\ee
Due to Lemma \ref{record}, the minimum value of $t\rightarrow Z_m(\lambda, t)$ on
$(-\infty,\pi_m]$ is  the value of $Z_m(\lambda,.)$ on the interval $\left[g_m(\lambda),
g_m\left(\lambda+\frac{1}{m}\right)\right)$. On the other hand,
due to the Cameron-Martin-Girsanov formula, almost surely,
$t\rightarrow z(t)$ has only one minimum on the interval $(-\infty,
\rho_1]$, but by definition of $\Psi_{\lambda}$ this unique minimum is at
 $g(\lambda)$.  Thus, by uniform  convergence of
 $Z_m(\lambda, t)$ to $z(t)$, almost surely, $\lim g_m(\lambda)=g(\lambda)$.
Still by definition of $\Psi_{\lambda}$,  $d(\lambda)$ is  the first hitting time of
level $z\left(g(\lambda)\right)$ after $\rho_1$:
\[d(\lambda)=\inf\{t\ge\rho_1\ :\ z(t)\ge z\left(g(\lambda)\right)\},\]
but due to   Proposition \ref{record},
\be
Z_m\left(\lambda, g_m(\lambda)\right)&=&Z_m\left(\lambda, d_m(\lambda)\right)
+\frac{1}{\sqrt{m-\lceil\lambda\sqrt m
\rceil}},\\
\lim_m Z_m\left(\lambda, d_m(\lambda)\right)&=&z\left(g(\lambda)\right),
\ee
so that $\liminf d_m(\lambda)\ge d(\lambda)$. Because  $d(\lambda)$ is  a stopping time,
almost surely there exists a sequence $t_k\downarrow d(\lambda)$ such that
\[z(t_k)<z(d(\lambda))=z(g(\lambda)).\]
 Thus
\be
\lim_m Z_m\left(\lambda, t_k\right)&<&\lim_m Z_m\left(\lambda, 
d_m(\lambda)\right),
\ee
and $d_m
(\lambda)<t_k$ for $m$ large enough.
Finally,  almost surely,
\begin{eqnarray}
\nonumber
\lim_m\ d_m(\lambda)&=&d(\lambda),\\
\label{bornes}
\lim_m\ R_1^{(m)}(\lambda)&=&R_1(\lambda).
\end{eqnarray}

\subsection{Proof of Theorem \ref{decomp}}
\label{p3}

We essentially do the same surgery on $h_m(t)=\psi_m(\lambda,t)$
as we did on   $\Psi_{\lambda}e$ at Figure 4: the analogs, for $h_m$,
of  properties $(i)$, $(ii)$ and $(iii)$ of Theorem \ref{decomp}
 are  combinatorial properties of parking schemes.
 Going to the limit  with the help
 of Theorem \ref{monticule} then yields
Theorem \ref{decomp}. We set
 \[n=\lfloor m-\lambda \sqrt m\rfloor\textrm{ and }
\ell=\lceil \lambda \sqrt m\rceil,\]
and we  consider   the probability space of  Subsection \ref{prfidis2}, enlarged
 to obtain a uniform random variable $w$, independent of
 $\left((\alpha_m)_{m\geq 1},b,e,v,(u_i)_{i\geq 1}\right)$.
From a parking scheme of $n$ cars on $m$ places, generated
with the help of  $\alpha_{m,\ell}$ as in Section \ref{uniform}, we obtain
 a profile  $h_m(t)=\psi_m(\lambda,t)$, and we have in mind to decompose it as
shown on Figure 9:  extend $h_m$ to a periodic function on the line, and set
\begin{eqnarray*}
q_m&=& \left(\tau_{v_m(\lambda)}h_m\right)^{[g_m(\lambda),d_m(\lambda)]},\\
r_m(t)&=& \left(\tau_{v_m(\lambda)}h_m\right)^{[d_m(\lambda),1+g_m(\lambda)]}\\
w(m)&=&\frac{\left\lceil \left(m-1-R^{m,\ell}_1\right)w\right\rceil}
{m-1-R^{m,\ell}_1}.
\end{eqnarray*}

\begin{figure}[ht]
\begin{center}
\includegraphics[width=9.7cm]{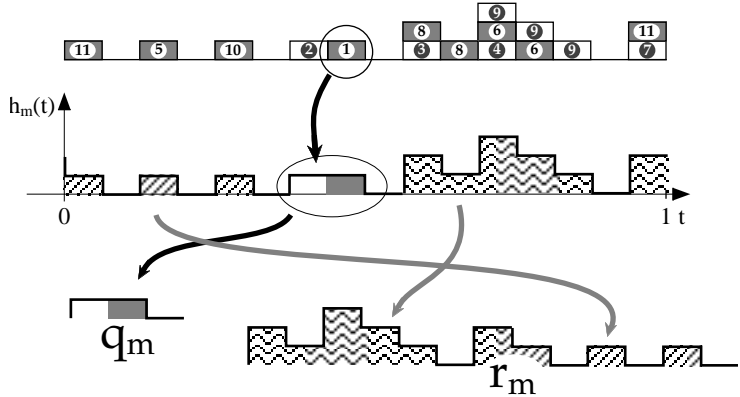}
\end{center}
\label{thegap2}
\caption{Decomposition  of the profile $h_m$}
\end{figure}

From  relation (\ref{combi})  ($n\leq m-2$), there are
\[C_{n-1}^{k-1}m(k+1)^{k-1}(m-k-1)^{n-k-1}(\ell-1)\]
parking schemes such that the block containing
 $c_1$ has $k$ cars ; for $C_{n-1}^{k-1}m$  of these parking schemes, the $k$ cars
of the block containing $c_1$ are parked according to a given parking scheme $\Pi$,
and  the remaining $n-k$ cars  are parked according to another given parking
 scheme $\tilde{\Pi}$: $C_{n-1}^{k-1}$ choices for the elements of the block containing
 $c_1$, $m$ choices for the position of this block. Thus, according to (\ref{combi}),
 the conditional probability of the parking schemes $(\Pi,\tilde{\Pi})$,
given that  $mR^{(m)}_1(\lambda)=k$,  is
\[\frac{1}{(k+1)^{k-1}}\times
\frac{1}{(m-k-1)^{n-k-1}(m-n-1)}.\]
That is, we have:
\begin{pro}
\label{discanalog}
Given   $R^{(m)}_1(\lambda)$, $r_m$  and $q_m$  are the profiles of
independent random uniform    confined parking schemes.
\end{pro}

This Proposition is the discrete analog of Theorem \ref{decomp},
so, to end the proof, we just have to go (carefully) to the limit.
In order to do that we need additional notations:
let $C$ be the space of continuous functions on
$[0,1]$, with the topology of uniform
convergence, and let $D$ be the space of
 cadlag functions on the same interval, embedded
with the Skorohod topology (see \cite[Ch. 3]{BILL}). The triplet of independent
random variables $(\Psi_{\lambda}e,\rho_1,w)$
with value in $C\times [0,1]^2$ defines the random variable
$(R_1(\lambda),q,\tau_wr)$ and its law $Q$, that is a probability measure
on the space $[0,1]\times C^2\subset [0,1]\times D^2$. The normalized
Brownian excursion
$e$ $\left(\textrm{resp. }\tau_w\left(\Psi_{\frac{\lambda}{\sqrt{1-x}}}e\right)\right)$
 defines the probability measure $\nu$
(resp. $\mu_x$)  on $C$. Now Theorem \ref{decomp}
is equivalent to:
\begin{eqnarray}
\label{aaa}
\int_{0}^1
f(\lambda,x)\int_{C^2}\Phi(x, y, z)
\mu_x(dz)\nu(dy)dx&=&E[\Phi(R_1(\lambda),q,\tau_wr)]\\
\nonumber
&=&\int_{R\times C^2}\Phi(x, y, z)Q(dx,dy,dz),
\end{eqnarray}
for any bounded uniformly continuous function $\Phi$
 on the space $[0,1]\times D^2$. It is harmless to assume
that $\Phi=0$ outside
$[a,b]\times D^2$, for some choice $0<a<b<1$.
On a probability space
$(\Omega,A,P)$, we already exhibited
a triplet  $(e,\rho_1,w)$
and  a sequence of $[0,1]\times D^2$-valued
random variables $\left(R^{(m)}_1(\lambda),q_m,\tau_{w(m)}r_m\right)$, satisfying:
\begin{itemize}
\item[($1$)]  almost
surely in $\Omega$, $\left(R^{(m)}_1(\lambda),q_m,\tau_{w(m)}r_m\right)$ converges
 to $\left(R_1(\lambda),q,\tau_wr\right)$  for the product topology
of  $[0,1]\times D^2$;
\item[($2$)] $\Pr\left(R^{(m)}_1(\lambda)=\frac{k}{m}\right)=
\varphi(m,n,k)$, in which
 $n=\lceil m-\lambda \sqrt m\rceil$;
\item[($3$)] the conditional law, $\nu_k$, of $q_m$  given that
 $R^{(m)}_1(\lambda)=\frac{k}{m}$, does not depend on $m$ and satisfies:
\[\nu_k\build{\tend}{}{weakly}\nu;\]
\item[($4$)] the conditional law, $\mu_{m,k}$, of $\tau_{w(m)}r_m$ given that
 $R^{(m)}_1(\lambda)=\frac{k}{m}$, satisfies:
\[\mu_{m,k}\build{\tend}{}{weakly}
\mu_x\hspace{1cm}\textrm{if }
m\tend \infty\hspace{0,1cm}\textrm{and }
k/m\tend x\in]0,1[;\]
\item[($5$)] conditionally, given that  $R^{(m)}_1(\lambda)=\frac{k}{m}$,
$\tau_{w(m)}r_m$  and $q_m$  are independent.
\end{itemize}

As $R^{(m)}_1(\lambda)=d_m(\lambda)-g_m(\lambda)$,
convergence of the first component  in point ($1$) is just
(\ref{bornes}). Uniform convergence of $h_m$ to $\tau_v\Psi_{\lambda}e$,
 uniform continuity of $\Psi_{\lambda}e(t)$ and (\ref{bornes})
 entails the uniform convergence of $q_m$ to $q$ and the uniform convergence
 of $\tau_{w(m)}r_m$ to $\tau_{w}r$, completing point $(1)$. Point ($2$) is just
  relation (\ref{combi}).
As a consequence  of Proposition \ref{discanalog}, given that $R^{(m)}_1(\lambda)=k$, $q_m$  is the
 profile of a random uniform confined  parking scheme of  $k$ cars on $k+1$ places,
so $q_m$ converges weakly to  a normalized Brownian excursion, as a special case
$\lambda=0$ of Theorem \ref{monticule} (see also \cite[Section 4]{CM}). That is,
$\nu_k$ converges weakly to $\nu$, giving point  $(3)$.
Similarly, given that $R^{(m)}_1(\lambda)=k$, $r_m$ is the profile of a  random uniform confined
 parking scheme for the $n-k$ remaining cars on the $m-k-1$ remaining places, and
$(m-k-1)w(m)$   is random uniform on $\{1,2,\ \dots\ ,m-k-1\}$, so $\tau_{w(m)}r_m$ is the profile
 of a  random  uniform \textit{non-confined}
 parking scheme of $n-k$  cars on  $m-k-1$  places.
If $k \simeq xm$, this parking scheme has \[\ell-1=m-n-1 \simeq \lambda \sqrt m
 \simeq \frac{\lambda}{\sqrt{1-x}}\ \sqrt{m-k-1}\] empty places,
thus Theorem \ref{monticule} applied to the  conditional law
$\mu_{m,k}$  of $\tau_{w(m)}r_m$ yields point  $(4)$.  Point  $(5)$
is already contained in Proposition \ref{discanalog}.

As a consequence of ($1$):
\begin{eqnarray}
\nonumber
\lim_mE\left[\Phi(R^{(m)}_1(\lambda),q_m,\tau_{w(m)}r_m)\right]&=&
E\left[\Phi(R_1(\lambda),q,\tau_wr)\right],
\end{eqnarray}
for any bounded uniformly continuous function $\Phi$.
We shall prove now that  properties ($2$) to ($5$)
 are sufficient to insure that, for any choice $0<a<b<1$, and for any bounded
 uniformly continuous function $\Phi$
satisfying $\Phi=0$ outside
$[a,b]\times C^2$,  we have
\begin{eqnarray}
\nonumber
\lim_mE[\Phi(R^{(m)}_1(\lambda),q_m,\tau_{w(m)}r_m)]&=&\int_{0}^1
f(\lambda,x)\int_{C^2}\Phi(x, y, z)
\mu_x(dz)\nu(dy)dx,
\end{eqnarray}
entailing (\ref{aaa}).

Let $M$ be a bound for $|\Phi|$ . Set:
\begin{eqnarray*}
E&=&\int_{0}^1f(\lambda,x)\int_C
\int_C\Phi(x, y, z)\mu_x(dz)\nu(dy)dx\\
&=&
E\left[\Phi(R_1(\lambda),q,\tau_wr)\right]\\
E_{1,m}&=&\int_{a}^bf(\lambda,x)\int_C
\int_C\Phi(x, y, z)\mu_{m,\lceil mx\rceil}
(dz)\nu_{\lceil mx\rceil}(dy)dx\\
E_{2,m}&=&\int_{a}^bf(\lambda,\lceil mx\rceil/m)\int_C
\int_C\Phi(\lceil mx\rceil/m, y, z)\mu_{n,\lceil mx\rceil}
(dz)\nu_{\lceil mx\rceil}(dy)dx\\
&=&\frac{1}{m}\sum_{k=1}^mf(\lambda,k/m)\int_C
\int_C\Phi(k/m, y, z)\mu_{m,k}(dz)\nu_k(dy)dx\\
E_{3,m}&=&\sum_{k=1}^m\varphi(m,n,k)\int_C
\int_C\Phi(k/m, y, z)\mu_{m,k}(dz)\nu_k(dy)dx\\
&=&E[\Phi(R^{(m)}_1(\lambda),q_m,\tau_{w(m)}r_m)].
\end{eqnarray*}

The last equality is a consequence of point  $(5)$.
By dominated convergence, owing to ($3$)
and  ($4$), $\lim_m E_{1,m} =E$.
By uniform continuity of $q$ and $\Phi$, $\lim_m E_{1,m}-E_{2,m} =0$.
Finally $\lim_m E_{2,m}-E_{3,m} =0$ due to Lemma \ref{local}.
$\hspace{1cm}\diamondsuit$

\section{Parking, fragmentation processes and
 the standard additive coalescent}
\label{addcoal}

In this Section, we give the proofs of   Theorems \ref{mixed},
 \ref{Levy} and  \ref{id}. These results are  consequences of Theorem
\ref{decomp}.

\subsection{Proof of Theorem \ref{mixed}}
\label{p1.5}

It should  be possible, following the line of proof
of  Subsection \ref{p3}, to exhibit a space
$(\Omega,A,P)$ on which there is almost sure convergence
of $\frac{1}{m}(R_1^m,R_2^m,\ ...\ ,R_k^m)$ to
$(R_1(\lambda),R_2(\lambda),\ ...\ ,R_k(\lambda))$
for each $k$, therefore yielding Theorem \ref{mixed}.
We rather borrow the clever  idea of
  \cite[Section 6.4]{PITYOR},
 that uses the decomposition of sample paths
 of a Brownian bridge
to compute the distribution
 of the sequence of widths of its excursions
(in that case a Poisson-Dirichlet distribution).

We introduce, as in \cite{PITYOR},
 a sequence $\rho=(\rho_k)_{k\geq 0}$ of
uniform random variables,  $\rho$ being
independent of $e$. With probability $1$, $\Psi_{\lambda}e(\rho_k)>0$:
if the excursion containing $\rho_k$  has width $B_j(\lambda)$,  we define
\[I_k=j,\]
yielding a size-biased permutation
of $B(\lambda)$, as explained in the introduction.
 Set:
\be
T(1)&=&\inf\left\{i\geq 2~|~\rho_i\notin [g(\lambda),d(\lambda)]\right\}\\
T(k+1)&=&\inf\left\{i\geq T(k)+1~|~\rho_i\notin [g(\lambda),d(\lambda)]\right\}.
\ee
The random variables $\rho_{T(k)}$ are independent and uniformly distributed
 on $[0,g(\lambda)]\cup[d(\lambda),1]$, and, almost surely, there exist
a unique number $\tilde{\rho}_k\in]0,1[$ such that
\[\rho_{T(k)}=\left\{d(\lambda)+\tilde{\rho}_k(1-R_1(\lambda))\right\} \hspace{0,1cm};\]
$\tilde{\rho}=(\tilde{\rho}_k)_{k\geq 1}$ is a sequence
of independent  random variables, uniform on $[0,1]$,
 and independent of $(e,\rho_1)$. Set  $\theta R(\lambda)=
\left(R_k(\lambda)\right)_{k\geq 2}$:
$(\rho_1,e)$ defines $R_1(\lambda)$,
 but among
the $(\rho_i)_{i\geq 2}$, only the $\rho_{T(k)}$
are useful  to determine $\theta R(\lambda)$.
Actually, up to  a multiplicative factor $1-R_1(\lambda)$,
 $\theta R(\lambda)$ is the size-biased permutation,
built with the help of the sequence $\tilde{\rho}$,
of the sequence of  widths of excursions of $r$,
or, equivalently, the size-biased permutation, built with the help of the sequence
$\hat{\rho}=\left(\{\tilde{\rho}_k-w\}\right)_{k\geq 1}$,
of the sequence of widths of excursions of $\tau_wr$.
Clearly $\hat{\rho}$ is a sequence of independent
and uniform random variables,
independent of $(r,w)$.
In view of Theorem \ref{decomp}($iii$), this leads to
\begin{lem}
\label{induc}
Given that
$R_1(\lambda)=x$, the sequence
 $\theta R(\lambda)=(R_k(\lambda))_{k\geq 2}$
 is distributed as $(1-x)R\left(\frac{\lambda}{\sqrt{1-x}}\right)$.
\end{lem}

Set:
\[s_k=x_1+\cdots+x_k.\]
Using  Lemma \ref{induc}, an easy
induction  on $k$ gives
the two following properties
\begin{itemize}
\item[(1)]
 $(R_j(\lambda))_{1\leq j\leq k}$ has
 the   distribution asserted in Theorem \ref{mixed} ;
\item[(2)] conditionally, given  that $(R_j(\lambda))_{1\leq j\leq k}
=(x_j)_{1\leq j\leq k}$,
$\theta^kR(\lambda)$
 is distributed as $(1-s_k)R\left(\frac{\lambda}{\sqrt{1-s_k}}\right)$,
\end{itemize}
ending the proof.
 $\hspace{1cm}\diamondsuit$

\subsection{Proof of Theorem \ref{Levy}}
\label{p2}

The operator   $\Psi_{\lambda}$ has the semigroup property, and, if $a$ and $a+x$ are
 two zeroes of a nonnegative function $f$, due to the Brownian scaling,
\[(\Psi_{\lambda}f)^{[a,a+x]}=\Psi_{\lambda\sqrt x}\left(f^{[a,a+x]}\right).
\]
Thus, conditionally, given that $\left(g(\lambda),d(\lambda)\right)=(a,a+x)$,  we have:
\be
\left(\Psi_{\lambda+\mu}e\right)^{[g(\lambda),d(\lambda)]}&=&
\left(\Psi_{\mu}(\Psi_{\lambda}e)\right)^{[g(\lambda),d(\lambda)]}\\
&=&\Psi_{\mu\sqrt{x}}\left((\Psi_{\lambda}e)^{[g(\lambda),d(\lambda)]}\right)\\
&=&\Psi_{\mu\sqrt{x}}\ q.
\ee
Still conditionally, $\rho_1$ is distributed
 as $a+x\rho_2$ in which $\rho_2$ is uniform on $]0,1[$, not
depending on $a$.
Thus, Theorem \ref{decomp} (ii), with  $(q,\mu\sqrt x,\rho_2)$  replacing
$(e,\lambda,\rho_1)$,
entails that, given $R_1(\lambda)=x$,
$(R_1(\lambda+\mu))_{\mu\geq 0}$
 is distributed as $
(x\hspace{0,1cm}R_1(\mu\sqrt x))_{\mu\geq 0}$.
Equivalently, by change of variables,
the conditional distribution
of $(\Sigma(\lambda+\mu))_{\mu\geq 0}$,
given that $\Sigma(\lambda)=y$,
is the same as the unconditional distribution
of \[\left((1+y)\Sigma\left(\frac{\mu}{\sqrt {1+y}}\right)+y\right)_{\mu\geq 0}.\]

This last statement yields (\ref{fidis}), by induction on $k$:
assuming that property at rank $k-1$ holds, we see that,
given that $\Sigma(\lambda_1)=y$,
\be
(\Sigma(\lambda_1+\lambda_2+ ... +\lambda_i))_{2\leq i \leq k}
&\build{=}{}{law}&\left((1+y)\Sigma
\left(\frac{\lambda_2+ ... +\lambda_i}{\sqrt {1+y}}\right)
+y\right)_{2\leq i \leq k}\\
&\build{=}{}{law}&\left(y+\frac{\lambda_2^2}{N_2^2}+
 ... +\frac{\lambda_i^2}{N_i^2}\right)_{2\leq i \leq k}.
\ee
Owing to Theorem \ref{decomp},
\[\Sigma(\lambda_1)
\build{=}{}{law}
\frac{\lambda_1^2}{N_1^2}.\hspace{0.5cm}\diamondsuit
\]

\subsection{Parking and the additive coalescent}
\label{p1}

In this subsection, we give an alternative proof of Bertoin's Theorem \ref{id}:
the coalescence of excursions of $\Psi_{\lambda}$,  as $\lambda\searrow 0$,
 has the same law  as  the coalescence of continuum random trees
 in the standard additive coalescent of Aldous and
Pitman \cite{ADD}. As we do not claim novelty,
our  proof will be sketchy at some points, but we hope to show that
some properties of the additive coalescent seem natural,
 once translated to parking schemes.

First let us prove that $B$  has the Markov property, that is:
\begin{equation}
\label{markov2}
\Pr\left(B(\lambda+\tilde\lambda)\in A\hspace{0.1cm}|\hspace{0.1cm}
{\mathcal G}_{\lambda}\right)=
\Pr\left(B(\lambda+\tilde\lambda)\in A\hspace{0.1cm}|\hspace{0.1cm}
B(\lambda)\right),
\end{equation}
in which ${\mathcal G}_{\lambda}$ is a sigma-field that contains
all the information about $\left(B(\lambda)\right)_{0\le\mu<\lambda}$ (the past of the
 process). Following closely \cite[Section 2]{BERT}, let ${\mathcal G}_{\lambda}$
stand for the P-completed sigma-field generated by
$\left(S^{\lambda}_t\right)_{0\le t\le 1}$, in which
\[S^{\lambda}_t=\sup_{0\le s\le t} \lambda s-e(s).\]
Bertoin \cite{BERT} argues that
 the complement of the support
of the Stieltjes measure $dS^{\lambda}$
is the union of nonoverlapping open intervals whose
lengths are given by $B(\lambda)$,
making $B(\lambda)$
 ${\mathcal G}_{\lambda}$-measurable. Bertoin gives furthermore
the following Skorohod-like formula, for $0\le\mu<\lambda$,
\[S^{\mu}_t=\sup_{0\le s\le t}S^{\lambda}_s-(\lambda-\mu)s,\]
with the consequence that $B(\mu)$
is ${\mathcal G}_{\lambda}$-measurable too.
Chassaing and Janson  \cite{CJ} give a construction of
$\Psi_{\lambda}e$ from the three independent sequences $U$, $\eta$ and $B(\lambda)$,
where  $\eta=(\eta_k)_{k\ge 1}$ is a sequence
 of independent Brownian excursions and $U=(U_k)_{k\ge 1}$  is a sequence
of independent  random variables uniform on $[0,1]$.
In this construction,
$S^{\lambda}$ depends only on $U$ and $B(\lambda)$, while the way
$B(\lambda)$ is fragmented to give $B(\lambda+\tilde\lambda)$ depends only  on
 $B(\lambda)$ and $\eta$, yielding  (\ref{markov2}).

We give some details, because the construction  of \cite{CJ}
describes the distribution of the point process $\Theta(\lambda)$ (see Section \ref{uniform})
that keeps track of  positions of excursions, and as such, this construction
gives some light on the distribution of positions of blocks: in the limit model,
the set of zeroes of $\Psi_{\lambda}e$ (empty places) is so to say Cantor-like, as there
 infinitely many excursions (blocks) between any given pair of zeroes.
First we build a copy $\left(X_{\lambda}(t)\right)_{0\le t\le 1}$
of the reflected Brownian bridge $|b|$ conditioned on its local
 time at zero $L_1(b)=\lambda$, following \cite[Section 6]{ASIN}: we place side
by side excursions with shape $\eta_i$ and width  $B_i(\lambda)$, the order of excursions
 being  dictated by $U$, that is, the excursion
with width $B_i(\lambda)$ and shape $\eta_i$ is on the left of the excursion
$(B_j(\lambda),\eta_j)$ if  $U_i<U_j$. To be formal, set
\be
g_i&=&\sum_{j, U_j<U_i}B_j(\lambda)\\
d_i&=&\sum_{j, U_j\le U_i}B_j(\lambda)\\
&=&g_i+B_i(\lambda),
\ee
 and let $X_{\lambda}$ be defined, on a dense subset of $[0,1]$, by
\[X_{\lambda}^{[g_i,d_i]}=\eta_i.\]
Then $X_{\lambda}$ is extended by continuity to $[0,1]$.
Incidentally, this random ordering of excursions is analog to the random insertion of
blocks in a confined parking scheme (cf. Subsection \ref{prsize-biased}): the analogy is used
in \cite{CJ} to prove that the limit of profiles of confined parking schemes is $X_{\lambda}$.
Then, in the same way as the random rotation of a random \textit{confined}
parking scheme gives a random parking scheme, $\Psi_{\lambda}e$
is obtained by random rotation of $X_{\lambda}$.

The proof of the Markov property requires  that this random rotation depends only on
$U$ and $B(\lambda)$, not on $\eta$, so let us define it.
According to \cite[Section 6]{ASIN}, the local time at $0$ of
$X_{\lambda}$, denoted $\left(L_t(X_{\lambda})\right)_{0\le t\le 1}$,
is defined for $t\in [g_i,d_i]$, by
\[L_t(X_{\lambda})=\lambda U_i,\]
thus $L_t$
depends only on $U$ and $B(\lambda)$. According to
\cite[Theorem 2.6 (i)]{CJ}, almost surely, there exists
a unique point $v$ in $[0,1)$ such that
\[L_v(X_{\lambda})-\lambda v=\max_{0\le t
\le1}
 L_t(X_{\lambda})
-\lambda t,\]
and $t\rightarrow X_{\lambda}(\{v+t\})$
is a copy of $\Psi_{\lambda}e$, that is, $\Psi_{\lambda}e$
is obtained from $X_{\lambda}$
through a random rotation $v$ that depends only
 on $U$ and $B(\lambda)$, not on $\eta$. Thus the local time at $0$ of
$\Psi_{\lambda}e$, $\left(L_t(\Psi_{\lambda}e)\right)_{0\le t\le 1}$,
is deduced from  $\left(L_t(X_{\lambda})\right)_{0\le t\le 1}$
through the shift $v$ as well, and depends only
 on $U$ and $B(\lambda)$.
Finally, according to
\cite[Proposition 8.2]{CJ},
\[L_t(\Psi_{\lambda}e)=S^{\lambda}_t,\]
so that ${\mathcal G}_{\lambda}$ is a subset of $\sigma(B(\lambda),U)$, the
P-completed sigma-field generated by $B(\lambda)$ and $U$.
As a by-product of the second part of the proof, we shall see that
 $$\sigma\left(B(\lambda+\tilde\lambda)\right)\subset\sigma(B(\lambda),\eta).$$
These two inclusions, with the independence between $B(\lambda)$, $\eta$ and $U$,
entail  the Markov property  (\ref{markov2}).

Once we know that both $B(\lambda)$ and  $Y(\lambda)$ have the Markov property,
we just have to check that they have the same transition probabilities, that is:
\begin{equation}
\label{markov3}
\Pr\left(B(\lambda+\tilde\lambda)\in A\hspace{0.1cm}|\hspace{0.1cm}
B(\lambda)=x\right)=
\Pr\left(Y(\lambda+\tilde\lambda)\in A\hspace{0.1cm}|\hspace{0.1cm}
Y(\lambda)=x\right).
\end{equation}
Let us describe the conditional distribution of $Y(\lambda+\tilde\lambda)$, given
$Y(\lambda)=x=(x_1,x_2,\,..)$. To this aim, let $\Delta_{s}$ denote
 the space of nondecreasing sequences $x$  of nonnegative numbers
with $\sum_{i\ge 1}x_i=s$ ; $B(\lambda)$ and  $Y(\lambda)$
are $\Delta_{1}$-valued random variables. Let $(y_k)_{k\ge 1}$ denote a
 sequence of independent random variables, $y_k$ being a $\Delta_{x_k}$-valued
 random variable with the same distribution as $x_k B\left(\tilde\lambda
\sqrt{x_k}\right)$ or as $x_k Y\left(\tilde\lambda
\sqrt{x_k}\right)$. According to  \cite[Section 1]{BERT},

\begin{pro}
\label{ptrsn}
Given $Y(\lambda)=x=(x_1,x_2,\cdots)$,  $Y(\lambda+\tilde\lambda)$ is distributed as the
decreasing rearrangement of the elements of sequences $y_1,y_2, \dots$\ \ .
\end{pro}
Let us prove that Proposition \ref{ptrsn} holds  also true with $B$ replacing $Y$, meaning,
informally,  that each of the clusters $x_i$ of the fragmentation process
$B(\lambda)$ starts anew a fragmentation process
distributed as $\left(x_iB(\tilde\lambda\sqrt x_i)\right)_{\lambda\geq 0}$.
We shall see that, in the case of $B$, the scaling factors $x_i$ and $\sqrt x_i$ come
from the Brownian scaling in the definition of $f^{[a,b]}$. These
  scaling factors can also be foreseen on  parking schemes:
 the time unit for the discrete fragmentation
process associated with parking $m$ cars on $m$ places,
is the departure of $\sqrt m$ cars. Due to the law of large numbers,
during one time unit, a given block of cars with size $x_im$
loses approximately $x_i\sqrt m=\sqrt x_i\sqrt {x_im}$ cars, meaning
 that,  for the internal clock of this block,
 $\sqrt x_i$ time units elapsed.

In order to give a formal proof, let $g_k$
(resp. $d_k=g_k+x_k$, $\eta_k$) denote the beginning (resp. the end, the shape) of
the excursion of $\Psi_{\lambda}e$ whose length is $B_k(\lambda)=x_k$:
\[\eta_k=\left(\Psi_{\lambda}e\right)^{[g_k,d_k]}.\]
As in Subsection \ref{p2}, we have:
\be
\left(\Psi_{\lambda+\tilde\lambda}e\right)^{[g_k,d_k]}
&=&\Psi_{\tilde\lambda\sqrt{x_k}}\left((\Psi_{\lambda}e)^{[g_k,d_k]}\right)\\
&=&\Psi_{\tilde\lambda\sqrt{x_k}}\ \eta_k.
\ee
Let $y_k$ denote the decreasing sequence of widths of those excursions of
$\Psi_{\lambda+\tilde\lambda}e$ that belong to the interval $[g_k,d_k]$:
$y_k$ is also, after normalisation by $x_k$,
the decreasing sequence of widths of excursions of
$\Psi_{\tilde\lambda\sqrt{x_k}}\ \eta_k$.
As consequences of Theorem \ref{decomp}, or of \cite{CJ,ASIN}, $B(\lambda)$
 and $\eta$ are independent,  and  $\eta$ is a sequence of
independent normalized Brownian excursions. Thus, given $B(\lambda)$,
 the sequences $y_k$ are independent
and respectively distributed as  $x_k B\left(\tilde\lambda
\sqrt{x_k}\right)$. We also have clearly
\[\sigma\left(B(\lambda+\tilde\lambda)
\right)\subset\sigma(B(\lambda),\eta).\]

\section{Concluding remarks }
\label{concl}

The additive coalescent has at least two constructions, seemingly quite different,
 one by Aldous and Pitman, through the continuum random tree,
the other by  Bertoin  through excursions of the family of stochastic processes $\left(\Psi_
{\lambda}e\right)_{\lambda\ge 0}$.
Actually Aldous and Pitman \cite{ADD} build the standard additive coalescent
 as the limit of a discrete model of coalescence-fragmentation:
they reverse the time of a discrete fragmentation
process  that starts with a random unrooted  labeled tree
 (the discrete analog of the continuum random tree),
whose edges are erased at random, one after the other  (the discrete analog of Poisson cuts).
In this paper we show that, asymptotically, parking schemes lead to Bertoin's construction
 of the additive coalescent.
As a first step towards a better understanding of the connection between these two different
 constructions of the additive coalescent, we show below an explicit connection between
 the discrete approximations for the additive coalescent, by random forests on one hand,
 and by parking schemes in the other hand.

  Erasing edges, Aldous and Pitman  \cite{ADD} split  the tree in smallest subtrees, and
 obtain a forest-valued stochastic process, but in \cite{ADD}  (as opposed to Pavlov's forests)
 the  forests are  \textit{unordered} sets of \textit{unrooted} trees. Then, focusing on the
 process of sizes of subtrees, and reversing time, Aldous and Pitman obtain a discrete
 Markovian coalescent, with  the following transition probability:
when the forest has $m$ nodes and $\ell$ subtrees the probability that two clusters
(subtrees) with sizes $x$ and $y$ merge in a larger
subtree with size $x+y$ is \cite[Lemma 1]{ADD}
\[\frac{x+y}{m(\ell-1)}\]
(as opposed to the fragmentation process where the edges are deleted
uniformly at random, in the time-reversed process, edges are not added
\textit{uniformly}).

There exists a striking similarity between the previous transition probability and relation
 (\ref{pt}), that gives the probability  of aggregation of two parking blocks with sizes $x$
 and $y$:  if we assume that the mass of a parking block with $x$ cars  is actually
$x+1$, for instance counting the empty place on the right of this block, then the process of
 sizes of blocks has the same distribution  as the discrete Markovian coalescent considered
in  \cite{ADD}.  We already knew that these two processes had the same \textit{asymptotic}
 distribution (the standard additive coalescent), but this is even better, and suggest the
 following question: are the two underlying richer structures,  the forest-valued stochastic
 process of on one hand, and the parking schemes-valued process on the other hand,
isomorphic in any sense ?

The answer is positive, up to a slight change in
Aldous \& Pitman' s model, that will be explained later. The description of the relation
 between the two discrete models
has several steps:  let    $\left(\Pi_k\right)_{0\le k\le m-2}$  denote the
\begin{figure}[ht]
\begin{center}
\includegraphics[width=14cm]{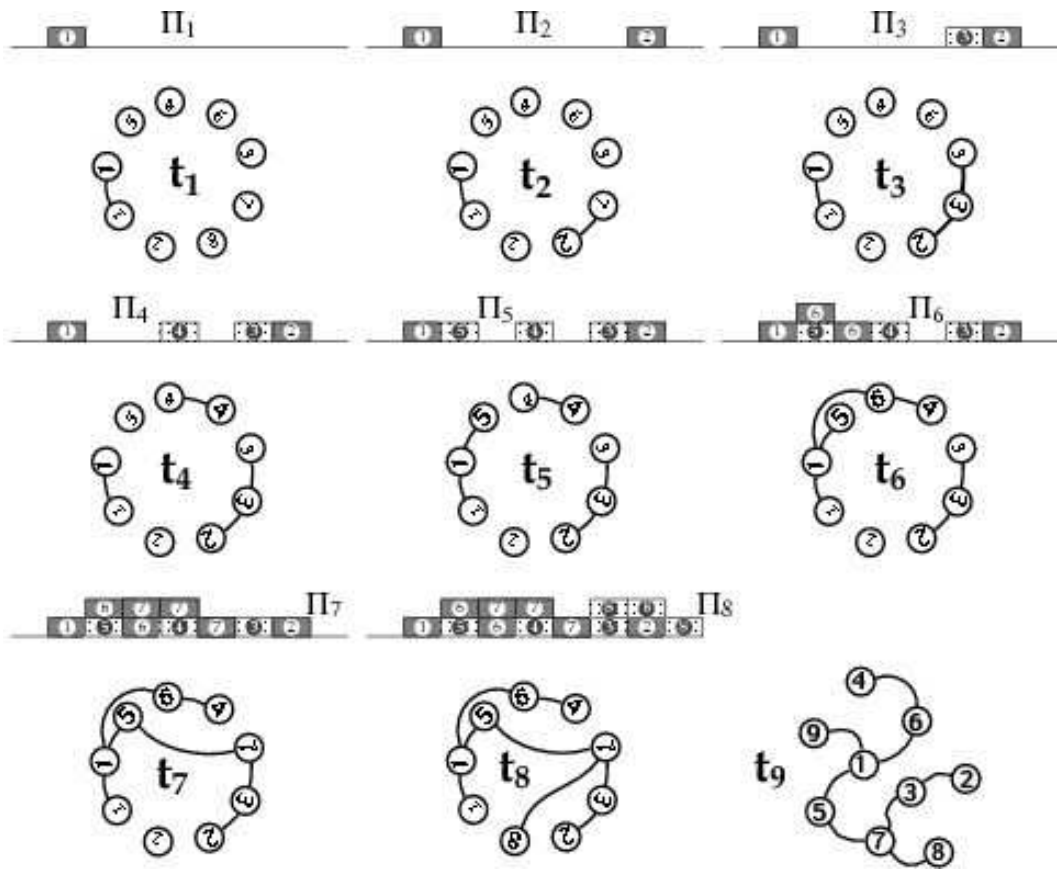}
\end{center}
\caption{Simultaneous fragmentation of trees and blocks.}
\end{figure}
 intermediate parking schemes leading to a confined parking scheme $\Pi_{m-1}$ for $m-1$
 cars on $m$ places,  through successive arrivals of cars. For $0\le k\le m-1$, let $t_k$
denote the Pavlov's forest
 associated with $\Pi_k$ through the one-to-one correspondence described at Subsection
 \ref{prPavlov}. Then $t_k$ can be obtained from $t_{m-1}$ by erasing edges in a
natural but \textit{deterministic} way. To describe it, consider  a random labeled tree $t_m$
 with $m$  vertices  $v_1,\ \dots\ ,v_m$. Relabel $v_m$ as $r_1$:  we get a Pavlov's tree
 $\tilde{t}_{m-1}$  with $m-1$ non-roots. Then erase the edge between $v_{m-1}$ and his
 father, and relabel $v_{m-1}$  as $r_2$: we get a Pavlov's tree $\tilde{t}_{m-2}$ with $m-2$
 non-roots. In  $\tilde{t}_{m-k}$, $v_{m-k}$ belongs to a subtree with root $r_i$, $i\le k$.
 To obtain  $\tilde{t}_{m-k-1}$  from $\tilde{t}_{m-k}$,  first,
 relabel roots $r_{i+1},\ \dots\ ,r_k$  as $r_{i+2},\ \dots\ ,r_{k+1}$,
 then erase the edge between $v_{m-k}$ and his father and relabel $v_{m-k}$ as $r_{i+1}$.
It is now easy to check that $\tilde{t}_{j}=t_j$.  In other terms, erasing at random edges of a
random rooted labeled tree, as in \cite{ADD}, or erasing successively the edge between
 $v_{m-1}$  (resp. $v_{m-2}$, $v_{m-3}$, ...) and its father, produces discrete Markovian
 coalescents with the same law, though the underlying forest-valued stochastic processes are
 different. We do not know if a different one-to-one correspondence forests-parking
would produce a nicer description of the forest-valued process associated to parking,
but very likely,
and as opposed to \cite{ADD}, given the random tree, the cuts would be deterministic.

In order to circumvent
this problem,  instead of drawing a random labeled tree, as in \cite{ADD}, we  draw
 the unlabeled \textit{random shape}  $st_m$ of a rooted labeled tree (that is,  $st_m$  is
a Galton-Watson tree with Poisson progeny, conditioned to have size $m$), then we delete
$st_m$'s  edges one by one
 in uniform random order. Compared with \cite{ADD}, the change is not fundamental.
But this last model is isomorphic to parking:  label  $v_{m}$ the root
of the random shape, then  label $v_{m-k}$ the vertice at the end (starting  from $v_m$)
 of the  $k^{th}$  deleted  edge, and we obtain a uniform random labeled tree $t_m$.
With him,
comes the associated random parking scheme $\Pi_{m-1}$. Focusing on sizes, and
eventually reversing
time, we get three  Markovian coalescents, from $st_m$, deleting edges at random,
from $t_{m-1}$, erasing successively the edge between
 $v_{m-1}$  (resp. $v_{m-2}$, $v_{m-3}$, ...) and its father, and reversing time,
or, finally, from $(\Pi_{k})_{0\le k\le m-1}$: these
three coalescents turn out to be equal, not only
in distribution, but also $\omega$ by $\omega$.

A (problematic) translation of this construction
to the continuous model  would open the way to a direct
 proof  of Theorem  \ref{id}, in which the Poisson
 process of cuts of the continuum random
tree, once  translated to Brownian excursion,
 would give the same cuts of excursions as those
 obtained from the  operators $\Psi_{\lambda}$.

Another possible direction for further research would be to explore connections between this
limit model for parking and the reflected Brownian motion with drift, or the Brownian storage
process, that appear as heavy traffic limits of queuing or storage systems \cite{HARR,lch1}.

\subsection*{Acknowledgements}
The starting point for this paper was the talk of Philippe Flajolet at
the meeting ALEA in February 1998 at Asnelles (concerning his paper
 with Viola \& Poblete), and a discussion that Philippe and the authors
 had in a small cafe of the French Riviera after the SMAI meeting of September
 1998. Some discussions of the first author with Marc Yor, and also with
 Uwe R\"osler, were quite fruitful. We thank  Philippe Lauren\c cot for calling
 our attention to the works of Aldous, Pitman \& Evans on coalescence models.
The papers of Perman, Pitman, \& Yor about random probability measures and
 Poisson-Dirichlet processes were also of a great  help. Finally  joint work
 \cite{CJ} with Svante Janson  lead to substantial changes  to a previous version
 of this paper. We are specially indebted to Svante for improvements to the proofs
of Subsection \ref{connect}, and he also pointed to us the connection with \cite{PAV}.
Finally we thank   R\'egine Marchand and two referees,
 whose careful reading led to substantial improvements.

\bibliography{resume}

\end{document}